\documentclass[11pt,a4paper,reqno]{amsart}
\usepackage{latexsym}
\usepackage{bbm,amssymb,amsthm,amsmath}
\usepackage{graphicx}
\usepackage{multirow}
\usepackage{caption}
\usepackage{bm}
\usepackage{multirow}
\usepackage{mathtools}  
\usepackage{float}  

\usepackage{xcolor}

\usepackage{tikz}	
\usetikzlibrary{shapes,arrows}
\usetikzlibrary{positioning}

\calclayout

\theoremstyle{plain}
\newtheorem{theorem}{Theorem}[section]	
\newtheorem{lemma}[theorem]{Lemma}
\newtheorem{corollary}[theorem]{Corollary}

\theoremstyle{definition}
\newtheorem{definition}[theorem]{Definition}
\newtheorem{example}[theorem]{Example}

\theoremstyle{remark}
\newtheorem{remark}[theorem]{Remark}

\numberwithin{equation}{section}	

\def\C{\mathbb{C}}
\def\R{\mathbb{R}}    
  
\def\N{\mathbb{N}}  

\def\Rd{\mathbb{R}^{d}}

\def\Cr{\mathbb{C}^{r}}
\def\Mat#1{{{\rm M}_{#1}(\C)}}

\let\mib=\boldsymbol

\def\mnu{{\mib \nu}}

\def\mx{\mathbf{x}}

\def\vf{\mathsf{f}}
\def\vu{\mathsf{u}}
\def\vv{\mathsf{v}}

\def\mA{\mathbf{A}}
\def\mB{\mathbf{B}}

\def\mI{\mathbf{I}}

\def\lD{\mathcal{D}}
\def\lH{\mathcal{H}}
\def\lV{\mathcal{V}}
\def\lW{\mathcal{W}}
\def\lX{\mathcal{X}}
\def\lY{\mathcal{Y}}
\def\lZ{\mathcal{Z}}

\def\phi{\varphi}

\def\supp{\operatorname{supp}}
\def\dom{\operatorname{dom}}
\def\dim{\operatorname{dim}}
\def\ker{\operatorname{ker}}
\def\ran{\operatorname{ran}}
\def\loc{\operatorname{loc}}
\def\cl{\operatorname{cl}}

\def\dup#1#2#3#4{{}_{#1\!}\langle\, #2 , #3 \,\rangle_{#4}} 
\def\scp#1#2{\langle\, #1 \mid #2 \,\rangle}  
\def\iscp#1#2{[\, #1 \mid #2 \,]}  

\begin{document}

\title[Classification of 1-d scalar Friedrichs operators]{Classification of classical Friedrichs differential operators: One-dimensional scalar case}

\author{M.~Erceg}\address{Marko Erceg,
	Department of Mathematics, Faculty of Science, University of Zagreb, Bijeni\v{c}ka cesta 30,
	10000 Zagreb, Croatia}\email{maerceg@math.hr}

\author{S.~K.~Soni}\address{Sandeep Kumar Soni,
	Department of Mathematics, Faculty of Science, University of Zagreb, Bijeni\v{c}ka cesta 30,
	10000 Zagreb, Croatia}\email{sandeep@math.hr}




\begin{abstract}
The theory of abstract Friedrichs operators, introduced by 
Ern, Guermond and Caplain (2007), proved to be a successful setting
for studying positive symmetric systems of first order partial 
differential equations (Friedrichs, 1958), 
nowadays better known as Friedrichs systems.
Recently, Antonić, Michelangeli and Erceg (2017) presented a purely
operator-theoretic description of abstract Friedrichs operators, 
allowing for application of the universal operator extension theory
(Grubb, 1968). 
In this paper we make a further theoretical step by developing a decomposition of the graph space
(maximal domain) as a direct sum of the minimal
domain and the kernels of corresponding adjoints.
We then study one-dimensional scalar (classical)
Friedrichs operators with variable coefficients
and present a complete
classification of
admissible boundary conditions.
\end{abstract}

\maketitle


\section{Motivation}

Friedrichs introduced
the concept of \emph{positive symmetric system} \cite{KOF}
(following his research on symmetric hyperbolic systems \cite{KOFh}), 
which are today customarily referred to as the \emph{Friedrichs system}.
More precisely, for a given open and bounded set  $\Omega \subseteq \Rd$ with Lipschitz boundary $\Gamma$, let the matrix functions 
$\mA_k \in W^{1,\infty}(\Omega;\Mat{r})$, $k=1,2,\dots,d$,
and $\mB \in L^\infty(\Omega; \Mat{r})$ satisfy 
\begin{equation}
\mA_k= \mA_k^\ast \qquad \hbox{on} \,\, \Omega
\tag{F1}
\end{equation}
and
\begin{equation}
(\exists{\mu_0>0}) \quad \mB + \mB^\ast+\sum_{k=1}^d 
    \partial_k \mA_k\geq 2\mu_0\mI
    \qquad {\hbox{a.e.~on}\,\, \Omega}\,.
\tag{F2}
\end{equation}
Then the first-order differential operator $L : L^2(\Omega)^r\longrightarrow \mathcal{D}'(\Omega)^r$
defined by
\begin{equation}
L \vu \;:=\; \sum_{k=1}^d \partial_k(\mA_k \vu) + \mB \vu
\tag{CFO}
\end{equation}
(here derivatives are taken in the distributional sense) 
is called \emph{the (classical) Friedrichs operator} or \emph{the symmetric positive operator}, while 
(for given $\vf \in L^2(\Omega)^r$) the first-order system of partial differential
equations $L \vu = \vf$ is called \emph{the (classical) Friedrichs system} or 
\emph{the symmetric positive system}.

Already Friedrichs showed that a wide variety of equations of mathematical physics
(regardless of their order), 
including classical elliptic, parabolic
and hyperbolic equations, can be adapted, or rewritten, 
in the form (CFO).
However, as pointed out explicitly by Friedrichs himself in \cite{KOF}, the main motivation of his 
approach `was not the desire for a unified treatment of elliptic and hyperbolic equations, but the 
desire to handle equations which are partly elliptic, partly hyperbolic, such as the Tricomi equation'.
One of the most important features of Friedrichs 
approach is in a clever way of representing different boundary (or initial) conditions by using
a matrix field on the boundary.
A nice historical 
exposition of the classical Friedrichs theory 
(which was very active until 1970's)
can be found in \cite{MJensen}.

New interest in Friedrichs systems arose from numerical analysis (see, for example, \cite{HMSW, MJensen}), thanks to their feature of providing a convenient unified
framework for numerical solutions to partial differential equations of different type.
However, well-posedness results obtained within the 
classical theory were not satisfactory;
there are only results on 
the existence of weak solutions, and the uniqueness of strong ones,
leaving the general question open on the joint existence and uniqueness of either a weak or a strong solution.
This motivated Ern, Guermond and Caplain \cite{EGC}
to introduce an abstract Hilbert space approach
for Friedrichs systems
(Definition \ref{def:abstractFO} below),
with an intrinsic formulation of the boundary 
conditions and a proper well-posedness result
(Theorem \ref{tm:well-posedness} below).
Thereafter this theory attracted the community 
for further theoretical 
and numerical investigations.
For example, studies of different representations of
boundary conditions and the relation with 
the classical theory 
\cite{ABcpde, ABjde, ABCE, ABVisrn, ABVnaRWA, AEM-2017, BH21},
applications to various (initial-)boundary value problems of elliptic, hyperbolic, and parabolic type
\cite{ABVjmaa, BH21, BEmjom, BVcpaa, DLS, EM19, EGsemel, MDS},
and the development of different numerical schemes \cite{BDG,TBT, BEF, EGsemel, EGbis, EGter}.

We are particularly interested in the results
obtained in \cite{AEM-2017}, where 
a purely operator-theoretic description of abstract 
Friedrichs operators is presented. 
In the new setting the authors proved that 
any abstract Friedrichs operator $T$ on a Hilbert space
$\lH$ admits a domain 
$\lV$ (i.e.~suitable boundary conditions), such that 
the abstract problem:
$$
\hbox{for a given $f\in\lH$ find $u\in\lV$ such that $Tu=f$}\;,
$$ 
is well-posed (this operator-theoretic reformulation of partial differential equations could be traced back to works of Vi\v{s}ik \cite{Vishik49, Vishik52}).
Moreover, it was recognised that Grubb's universal operator 
extension theory for the non-symmetric setting \cite{Grubb-1968}
is applicable (which could be seen as an extension and improvement 
of the results from \cite{Vishik52}), 
allowing for a complete classification of 
realisations of interest.
Our final aim is to further develop these results and 
then apply them to (CFO) (understood as an abstract Friedrichs 
operator) and provide a (possible complete) classification 
of all admissible boundary conditions.
As a first step, in this paper we develop 
a precise description of the graph space 
(maximal domain) in terms of a decomposition as a direct sum
of the minimal domain and the kernels of corresponding adjoints. 
This result gives better control over the choice of boundary conditions. Then we apply the theory to one-dimensional ($d=1$)
scalar ($r=1$) (CFO) with variable coefficients.

The paper is organised as follows.
In Section \ref{sec:abstractFO} we recall the definition and main properties of 
abstract Friedrichs systems, with the emphasis on recent insights.
The first important result of the paper is developed in Section 
\ref{sec:decomposition}, 
where we obtain a certain decomposition of the graph space
for general abstract Friedrichs operators
(Theorem \ref{tm:decomposition}). 
As a consequence, an explicit bijective realisation, with some nice additional
features, is recognised (Corollary \ref{cor:decomposition}).
In sections \ref{sec:1d-intro}--\ref{sec:examples}
we study one-dimensional ($d=1$) scalar  ($r=1$) classical Friedrichs differential operators
(CFO), where variable coefficients are allowed.
More precisely, in Section \ref{sec:1d-intro} we develop some 
preliminary results related to the graph space and the boundary operator.
In Section \ref{sec:1d-class} we give a complete classification of
admissible boundary conditions, i.e.~of bijective realisations,
based on Grubb's universal classification theory (briefly presented in Appendix).
The paper is closed by a few examples, presented in Section \ref{sec:examples},
illustrating the results of previous sections.
\smallskip

\noindent\textbf{Notation.} Most of our notations are standard, 
let us only emphasise the following. 
For the sake of generality, in the paper we work on 
complex vector spaces. Thus, by $\lH$ we denote 
a complex Hilbert space with scalar product 
$\scp\cdot\cdot$, which we take to be 
linear in the first and anti-linear in the second entry.
The corresponding norm is given by
$\|\cdot\|:=\sqrt{\scp\cdot\cdot}$.
For $\lH=\Cr$ we shall often use an alternative notation:
$\scp{\mathsf{x}}{\mathsf{y}} = \mathsf{x}\cdot\mathsf{y}$,
$\mathsf{x}, \mathsf{y}\in\Cr$.
The topological (anti)dual 
$\lH'$ will be identified with $\lH$ by means of the usual duality 
(the Riesz representation theorem).
For any Banach space $\lX$ by $\dup{\lX'}\cdot\cdot \lX$ 
we denote the corresponding dual product between 
$\lX$ and its (anti)dual $\lX'$.
The annihilator of $S\subseteq \lX$, denoted by $S^0$,
is a closed subspace of $\lX'$ given by 
$S^0=\{f\in\lX': (\forall u\in S) \ \dup{\lX'}{f}{u}\lX = 0\}$.
For a subspace $\lY\subseteq\lX$ we denote by 
$\cl_\lX \lY$ its closure within $\lX$.

For a densely defined linear operator $A:\lH\to \lH$ we denote by 
$\dom A$, $\ker A$, $\ran A$, $\overline{A}$, $A^\ast$ 
its \emph{domain}, \emph{kernel}, \emph{range} (or \emph{image}), 
\emph{closure} (if it exists), and \emph{adjoint}, respectively.
For $S\subseteq \lH$, the \emph{restriction} of $A$ to $S$ is denoted 
by $A|_S$.
For two linear operators $A, B$ in $\lH$ by $A\subseteq B$
we mean that $\dom A\subseteq \dom B$ and $B|_{\dom A}=A$.
By $\scp{\cdot}{\cdot}_A := \scp\cdot\cdot
+ \scp{A\,\cdot}{A\,\cdot}$ we denote the 
\emph{graph scalar product}, while the corresponding norm $\|\cdot\|_A:=\sqrt{\scp{\cdot}{\cdot}_A}$ 
is called the \emph{graph norm}. 
If $A=A^\ast$, 
then $A$ is said to be \emph{self-adjoint}, while the infimum of its spectrum is called the \emph{bottom}. 
The \emph{identity} operator is denoted by $\mathbbm{1}$.
For a \emph{direct} sum between two vector spaces we use the symbol $\dotplus$. 
We write $\ominus$ 
for the \emph{orthogonal difference} in order to express in which 
Hilbert space the orthogonal complement 
is taken. 

For any complex number $z\in\C$ we denote by $\Re z$ and $\Im z$
the real and the imaginary part of $z$, respectively.

\section{Abstract Friedrichs operators}\label{sec:abstractFO}

The idea behind the introduction of the abstract formalism 
is to express all important (pointwise) features 
of classical Friedrichs operators (CFO)
in an abstract setting.
Thus, the definition of abstract Friedrichs operators
should be broad enough to encompass classical Friedrichs
operators. However, not too broad, so that the required results can 
be obtained within the class, 
such as well-posedness.

The abstract Hilbert space formalism for Friedrichs systems which we 
study in this paper was introduced 
and developed in \cite{EGC, ABcpde} for real vector spaces, while 
the required differences for complex vector spaces have been 
supplemented more recently in \cite{ABCE}. 
Here we present the definition in the form given in
\cite[Definition 1]{AEM-2017}.

\begin{definition}\label{def:abstractFO}
A (densely defined) linear operator $T$ on a complex Hilbert space $\lH$ 
is called an \emph{abstract Friedrichs operator} if it admits another
(densely defined) linear operator $\widetilde{T}$ on $\lH$ with the following properties:
\begin{itemize}
 \item[(T1)] $T$ and $\widetilde{T}$ have a common domain $\lD$, 
 which is dense in $\lH$, satisfying
 \[
 \scp{T\phi}\psi \;=\; \scp\phi{\widetilde T\psi} \;, \qquad \phi,\psi\in\mathcal{D} \,;
 \]
 \item[(T2)] there is a constant $c>0$ for which
 \[
 \|(T+\widetilde{T})\phi\| \;\leqslant\; c\|\phi\| \;, \qquad \phi\in\mathcal{D} \,;
 \]
 \item[(T3)] there exists a constant $\mu_0>0$ such that
 \[
 \scp{(T+\widetilde{T})\phi}\phi \;\geqslant\; 2\mu_0 \|\phi\|^2 \;, \qquad \phi\in\mathcal{D} \,.
 \]
\end{itemize}
The pair $(T,\widetilde{T})$ is referred to as a \emph{joint pair of abstract Friedrichs operators}
(the definition is indeed symmetric in $T$ and $\widetilde{T}$).
\end{definition}

\begin{remark}
Another interesting abstract approach which covers similar differential operators
can be found in \cite{Pic09, PTW15}. This theory deals with an abstract operator
(instead of a pair), and in particularly it covers operators of the form
$\partial_t M_0 + M_1 + A$, where $M_0, M_1$ are bounded linear operators on $\lH$
and $A$ is an unbounded skew-self-adjoint operator on $\lH$ 
(e.g.~a first order differential operator in spatial variables).
For example, if (CFO) is a time-independent operator,
then one could take $M_0=0$, $M_1=\mB$ and 
$A = \sum_{k=1}^d \partial_k(\mA_k \,\cdot)$.

The main strength of this approach is in studying evolution problems,
but a small drawback is that operator $A$ should be independent of the 
time variable $t$, 
i.e.~coefficients of the differential operator in spatial variables 
should not depend on the time variable. 
A similar situation occurs in the non-stationary theory for 
abstract Friedrichs systems \cite{BEmjom}.
\end{remark}

The following characterisation of joint pairs of abstract Friedrichs
operators can be found in \cite[Theorem 8]{AEM-2017}.

\begin{theorem}
A pair of operators $(T,\widetilde T)$ on a complex Hilbert space $\lH$ 
is a joint pair of abstract
Friedrichs operators on $\lH$ if and only if
$T\subseteq \widetilde T^*$, $\widetilde T\subseteq T^*$, and
$\overline{T+\widetilde T}$ is an everywhere defined, bounded, self-adjoint operator on $\lH$ with strictly positive bottom. 
\end{theorem}

\begin{remark}
Condition (T3) is used in the previous theorem only to get that 
$\overline{T+\widetilde T}$ has strictly positive bottom. 
More precisely, a pair $(T,\widetilde T)$ satisfies 
conditions (T1)--(T2) if and only if  
$T\subseteq \widetilde T^*$, $\widetilde T\subseteq T^*$, and
$\overline{T+\widetilde T}$ is an everywhere defined, bounded, self-adjoint operator on $\lH$. 
Since many statements hold even in this case, we shall explicitly 
emphasise in which particular situations condition (T3) is necessary.
\end{remark}

Operators $A, B$ on $\lH$ with the property that $A\subseteq B^*$
and $B\subseteq A^*$ are often referred to as \emph{dual pairs}.
Thus, by the previous theorem, operators forming a joint pair of 
abstract Friedrichs operators are dual pairs (in fact this follows merely
from condition (T1)). 

Let $(T,\widetilde T)$ be a joint pair of abstract Friedrichs operators.
By (T1) it is evident that $T$ and $\widetilde T$ are closable. 
Since $T+\widetilde T$ is a bounded operator, graph norms 
$\|\cdot\|_{T}$ and $\|\cdot\|_{\widetilde T}$ are equivalent. 
The consequence is that (see \cite[Subsection 2.1]{EGC}
and \cite[Theorem 7]{AEM-2017})
\begin{equation}\label{eq:ww0}
\begin{aligned}
    \dom \overline T &\;=\; \dom \overline{\widetilde{T}} \;=:\; \lW_0 \;, \\
    \dom T^* &\;=\; \dom\widetilde{T}^* \;=:\; \lW \;,
\end{aligned}
\end{equation}
and $\bigl(\overline{T+\widetilde T}\bigr)|_\lW=\widetilde T^*+T^*$.
This implies that $(\overline T,\overline{\widetilde T})$ is also 
a pair of abstract Friedrichs operators.
Now we simplify our notation by introducing
\begin{equation*}
    T_0 \;:=\; \overline{T} \;, \quad 
    \widetilde{T}_0 \;:=\; \overline{\widetilde{T}} \;, \quad
    T_1 \;:=\; \widetilde{T}^* \;, \quad
    \widetilde{T}_1 \;:=\; T^* \;.
\end{equation*}
Therefore, we have
\begin{equation}\label{eq:dual-pairs}
    T_0 \;\subseteq\; T_1 \quad \hbox{and} \quad
    \widetilde{T}_0 \;\subseteq\; \widetilde{T}_1 \;.
\end{equation}

When equipped with the graph norm (one of two equivalent 
norms $\|\cdot\|_{T_1}$ and $ \|\cdot\|_{\widetilde{T}_1}$), 
the space $\lW$ becomes a Banach space, thus we shall call
it the \emph{graph space} (it is in fact a Hilbert space when we consider one of the graph scalar products $\scp{\cdot}{\cdot}_{T_1} $ or $\scp{\cdot}{\cdot}_{\widetilde{T}_1}$; however another inner product plays 
a more important role -- see (2.3) below).
On the other hand, $\lW_0$ is a closed 
subspace of the graph space $\lW$, 
while it is dense in $\lH$ (since it contains $\lD$).
As an illustration, for $\lH=L^2(\Omega)$ and 
a certain choice of operators we can achieve that $\lW$ and $\lW_0$
are Sobolev spaces $H^1(\Omega)$ and $H^1_0(\Omega)$, respectively.

In correspondence to the theory of symmetric operators
and having in mind applications to partial differential equations, 
the space $\lW$ could be called the \emph{maximal domain} 
(with no (initial-)boundary conditions prescribed), 
while $\lW_0$ the \emph{minimal domain} (with 
zero (initial-)boundary conditions).
This can be justified by using the boundary operator
(see Lemma \ref{lm:boundary_op}(ii) below), 
which serves also as the most natural way to study different boundary
conditions in this abstract setting.
The continuous linear map
\begin{align}
    & D : (\lW,\|\cdot\|_{T_1}) \to (\lW,\|\cdot\|_{T_1})' \nonumber \\
    \iscp uv \;:=\; & \dup{\lW'}{Du}{v}\lW \;:=\; \scp{T_1u}{v} 
        - \scp{u}{\widetilde{T}_1v} \;,
        \quad u,v\in\lW \;, \label{eq:D}
\end{align}
we call the boundary operator associated with the pair $(T_0,\widetilde{T}_0)$
(or equivalently with the pair $(T,\widetilde{T})$).
Some nice properties of the boundary operator we collect in the following lemma
(see \cite[Subsection 2.2]{EGC} and \cite[Lemma 1]{ABCE}).

\begin{lemma}\label{lm:boundary_op}
Let a pair of operators $(T,\widetilde{T})$ on $\lH$ satisfy 
{\rm (T1)--(T2)}. Then the boundary operator $D$ satisfies
 \begin{itemize}
    \item[i)] $(\forall u,v\in\lW) \quad
        \dup{\lW'}{Du}{v}\lW = \overline{\dup{\lW'}{Dv}{u}\lW} \,,$
    \item[ii)] $\ker D=\lW_0 \,,$
    \item[iii)] $\ran D=\lW_0^0 \,,$
 \end{itemize}
 where ${}^0$ stands for the annihilator.
\end{lemma}

The previous lemma ensures that $(\lW,\iscp \cdot\cdot)$
is an indefinite inner product space (see e.g.~\cite{Bo}).
Then $\iscp\cdot\cdot$-orthogonal complement of set $S\subseteq \lW$
is defined by
\begin{equation}\label{eq:orth-compl-D}
    S^{[\perp]} := \bigl\{u\in \lW : (\forall v\in S) \quad 
        \iscp uv = 0\bigr\}  \;,
\end{equation}
which is by definition a subspace of $\lW$.
Moreover, since $D$ is continuous, it is a closed subspace of $\lW$
with respect to the graph norm. 
For $L\subseteq S \subseteq \lW$ we have
$S^{[\perp]}\subseteq L^{[\perp]}$, while by 
Lemma \ref{lm:boundary_op}(i),(ii) it holds that $\lW_0^{[\perp]}=\lW$
and $\lW^{[\perp]}=\lW_0$
(this implies that $(\lW,\iscp \cdot\cdot)$ is a degenerate space).
Using the fact that the quotient space $\lW/\lW_0$ 
is a non-degenerate inner product space (more precisely 
a \emph{Kre\v\i n space}; \cite[Lemma 8]{ABcpde}) 
we have the following (see also Theorem IV and lemmas 7 and 9 
in the aforementioned reference and references therein). 
\begin{lemma}\label{lm:closed-compl}
A subspace $\lV$ of $\lW$ which contains $\lW_0$
is closed in $\lW$ (i.e.~with respect to the graph norm)
if and only if $\lV=\lV^{[\perp][\perp]}$.
\end{lemma}

\smallskip

For any subspace $\lV$ between the minimal and the maximal domain, 
i.e.~$\lW_0\subseteq \lV\subseteq \lW$, we call the restriction 
$T_1|_\lV$ a \emph{realisation} of $T_0$ (or $T$), and analogously
for $\widetilde T_1$. 
In terms of applications to partial differential equations, 
this can be seen as that each realisation corresponds to a different 
set of boundary conditions, which are prescribed implicitly by the 
choice of the domain $\lV$. 

Our main goal is to classify \emph{all} 
such closed subspaces $\lV$ (in $\lW$) for which we have that 
for any $f\in\lH$
the abstract problem $(T_1|_\lV) u=f$ is well-posed.
This implies that $T_1|_\lV$ is a closed densely defined 
\emph{bijective} 
operator on $\lH$. By the closedness, 
it is a $(\lV,\|\cdot\|_{T_1})\to \lH$ continuous map.
Thus, the inverse is a $\lH\to (\lV,\|\cdot\|_{T_1})$
continuous map, and hence continuous on $\lH$ as well
(see e.g.~\cite[Remark 6]{AEM-2017}).
An interesting geometrical consequence is that the adjoint
$(T_1|_\lV)^*$ has the same property. Indeed,
let us first note that from 
\begin{equation*}\label{eq:realisations}
T_0 \;\subseteq\; T_1|_\lV \;\subseteq\; T_1
\end{equation*}
we have
$$
\widetilde T_0 \;\subseteq\; (T_1|_\lV)^* \;\subseteq\; \widetilde T_1 \;.
$$
Thus, $\widetilde\lV:=\dom (T_1|_\lV)^*$ is a closed subspace in $\lW$
that contains $\lW_0$, and $(T_1|_\lV)^*=\widetilde T_1|_{\widetilde \lV}$.
By the standard results (cf.~\cite[Theorem 12.7]{GG}) 
operator $\widetilde{T}_1|_{\widetilde \lV}$ is also injective and has a range 
dense in $\lH$. Moreover, the range of $\widetilde{T}_1|_{\widetilde \lV}$ is
in fact closed since 
$$
(\widetilde{T}_1|_{\widetilde \lV})^{-1} = \bigl((T_1|_\lV)^*\bigr)^{-1}
    = \bigl((T_1|_\lV)^{-1}\bigr)^{*}
$$
is bounded on $\lH$ (the adjoint of a bounded operator is bounded).

The conclusion is the following:
\emph{if\/ $T_1|_\lV$ is a closed bijective realisation of $T_0$, 
then\/ $(T_1|_\lV)^*=\widetilde T_1|_{\widetilde \lV}$ is a
closed bijective realisation of\/ $\widetilde T_0$.}
Therefore, without any loss of generality we can simultaneously 
study both the original problem $T_1 u=f$ and the associated 
adjoint problem $\widetilde T_1 v = g$.
This means that our main goal can be reformulated:
we seek for pairs $(T_1|_\lV, \widetilde{T}_1|_{\widetilde \lV})$ 
of mutually adjoint \emph{bijective} realisations 
relative to $(T_0,\widetilde{T}_0)$.

In \cite[Lemma 11]{AEM-2017} the following useful 
characterisation in terms of the boundary operator
of mutual adjointness was derived. 

\begin{lemma}\label{lm:mutual-adj}
Let $(T_0,\widetilde T_0)$ be a pair of closed operators on 
$\lH$ satisfying conditions
{\rm (T1)--(T2)}, and let $\iscp\cdot\cdot$ be the 
associated indefinite inner product on the graph space $\lW$
given by \eqref{eq:D}.
Let $(T_1|_\lV,\widetilde T_1|_{\widetilde \lV})$ be a pair 
of realisations relative to $(T_0, \widetilde T_0)$,
i.e.~$\lV$ and $\widetilde \lV$ are subspaces of $\lW$
that contain $\lW_0$ (see \eqref{eq:ww0}). Then
\begin{equation*}
(T_1|_\lV)^* \;=\; \widetilde T_1|_{\widetilde \lV}
    \;\iff\; \widetilde \lV \;=\; \lV^{[\perp]}
\end{equation*}
and
\begin{equation*}
(\widetilde T_1|_{\widetilde \lV})^* \;=\; T_1|_{\lV}
    \;\iff\; \lV \;=\; \widetilde\lV^{[\perp]} \;.
\end{equation*}
In particular, if $\lV$ is closed in $\lW$, then 
condition $\widetilde \lV = \lV^{[\perp]}$ is sufficient to 
have that operators $T_1|_\lV$ and 
$\widetilde T_1|_{\widetilde \lV}$ are mutually adjoint.
\end{lemma}

Existence of such bijective mutually adjoint realisations 
was obtained in \cite[Theorem 13]{AEM-2017} and the result is 
the following.

\begin{theorem}\label{tm:existence-ref-op}
Let $(T_0,\widetilde T_0)$ be a joint pair of closed abstract Friedrichs
operators on $\lH$. 
\begin{itemize}
    \item[i)] There exists a pair 
    $(T_1|_\lV,\widetilde{T}_1|_{\widetilde{\lV}})$ of mutually adjoint
    bijective realisations relative to $(T_0,\widetilde{T}_0)$.
    \item[ii)] If both $\ker T_1\neq \{0\}$ and 
    $\ker\widetilde{T}_1\neq \{0\}$, then the pair
    $(T_0,\widetilde{T}_0)$ admits uncountably many mutually adjoint
    pairs of bijective realisations relative to $(T_0,\widetilde{T}_0)$.
    On the other hand, if either $\ker T_1=\{0\}$ or
    $\ker \widetilde{T}_1=\{0\}$, then there is exactly one mutually 
    adjoint pair of bijective realisations relative to
    $(T_0,\widetilde{T}_0)$. Such a pair is precisely
    $(T_1,\widetilde T_0)$ when $\ker T_1=\{0\}$, 
    and $(T_0,\widetilde T_1)$ when $\ker \widetilde{T}_1=\{0\}$.
\end{itemize}
\end{theorem}

The existence part of the previous theorem can be observed in the following
way:
any differential operator that can be cast into this theory of abstract 
Friedrichs operators admits at least one set of (initial-)boundary 
conditions for which the corresponding problem is well-posed.  
Existence of a pair of bijective realisations allows also for 
application of Grubb's universal classification theory, which was 
recognised in \cite{AEM-2017}, and which will be used in this 
manuscript as well.
Another consequence of the previous theorem is that $\ran T_1 = \ran \widetilde T_1 = \lH$. 

One can notice that in Theorem \ref{tm:existence-ref-op} 
condition (T3) is finally assumed,
since coercivity is important for the derivation of the result.
Indeed, the existence is established by using the following
sufficient condition which can be found already in the first paper
on abstract Friedrichs operators \cite[Theorem 3.1]{EGC}
(see \cite{ABCE} for the result in the complex setting).

\begin{theorem}\label{tm:well-posedness}
Let $(T_0,\widetilde{T}_0)$ be a joint pair of closed 
abstract Friedrichs operators on $\lH$,
and let $(T_1|_\lV,\widetilde{T}_1|_{\widetilde{\lV}})$ be
a pair of mutually adjoint realisations relative to 
$(T_0,\widetilde{T}_0)$.

If a pair $(\lV, \widetilde\lV)$ of subspaces of $\lW$ has 
\emph{definite sign} with respect to $\iscp\cdot\cdot$,
i.e.
\begin{equation}
\begin{aligned}
&(\forall u\in \lV) \qquad \iscp uu \;\geq\; 0 \,, \\
&(\forall v\in \widetilde\lV) \qquad \iscp{v}{v} \;\leq\; 0\,,
\end{aligned}
\tag{V1}
\end{equation}
where $\iscp\cdot\cdot$ is given in \eqref{eq:D},
then $(T_1|_\lV,\widetilde{T}_1|_{\widetilde{\lV}})$
is a pair of mutually adjoint bijective realisations relative to 
$(T_0,\widetilde{T}_0)$.
\end{theorem}

Coercivity is needed to apply the Banach-Nečas-Babuška
theorem (see \cite[Theorem 2.6]{EGbook}),
which ensures the required bijectivity.
With aid of condition (T3) we have the following result
(cf.~\cite[Lemma 3.2]{EGC}, \cite[Lemma 2]{ABCE}).

\begin{lemma}\label{lm:coercivity}
Let $(T_0,\widetilde{T}_0)$ be a joint pair of closed 
abstract Friedrichs operators on $\lH$
and let $(\lV,\widetilde{\lV})$ be a pair of linear subspaces of
the graph space $\lW$ satisfying condition {\rm (V1)}. 
Then operators $T_1|_\lV$ and 
$\widetilde{T}_1|_{\widetilde \lV}$
are $\lH$-coercive on $\lV$ and $\widetilde{\lV}$, respectively,
i.e.
\begin{align*}
& (\forall u\in\lV) \quad |\scp{T_1 u}{u}| \geq \mu_0 \|u\|^2 \,, \\
& (\forall v\in\widetilde\lV) \quad |\scp{\widetilde{T}_1 v}{v}| 
    \geq \mu_0 \|v\|^2 \,.
\end{align*}
\end{lemma}

\begin{remark}\label{ex:V1-kernels}
Let us consider a joint pair of closed abstract Friedrichs operators 
$(T_0,\widetilde{T}_0)$ on $\lH$.
\begin{itemize}
    \item[i)] A trivial pair satisfying condition (V1) is 
    $(\lW_0,\lW_0)$ since $\ker D=\lW_0$.
    This implies that closed operators 
    $T_0=T_1|_{\lW_0}$ and $\widetilde{T}_0=\widetilde{T}_1|_{\lW_0}$ 
    are $\lH$-coercive, hence injective.
    In particular, their ranges $\ran T_0$ and $\ran \widetilde T_0$
    are closed in $\lH$.

    Therefore, the following orthogonal decompositions of $\lH$ hold:
    \begin{equation}\label{eq:orth-decomp-lH}
    \begin{aligned}
    \mathcal{H} &\;=\;\ran T_0\oplus \ker \widetilde{T}_1 \\
    & \;=\; \ran \widetilde T_0\oplus \ker {T}_1 \;.
    \end{aligned}
    \end{equation}
\item[ii)] Here we present that the pair of subspaces 
$(\lW_0+\ker \widetilde T_1,\lW_0+\ker T_1)$ satisfies condition (V1).

Let $u_0\in\lW_0$ and $\tilde{\nu}\in\ker \widetilde{T}_1$ be arbitrary. 
We have
$$
\iscp{u_0+\tilde\nu}{u_0+\tilde\nu} = \dup{\lW'}{D\tilde\nu}{\tilde\nu}\lW
    = \scp{T_1\tilde{\nu}}{\tilde{\nu}}
    = \scp{(T_1+\widetilde{T}_1)\tilde{\nu}}{\tilde{\nu}}
    \geq 2\mu_0\|\tilde{\nu}\|^2 \geq 0 \,,
$$
where in the first equality we used Lemma \ref{lm:boundary_op}(i),(ii),
in the second and the third that $\widetilde T_1\tilde{\nu}=0$,
while the argument is closed by an application of the coercivity 
of $T_1+\widetilde T_1$ (i.e.~condition (T3)).
For the second subspace the calculation is completely analogous. 

Thus, the previous lemma implies that operators 
$T_1|_{\lW_0+\ker\widetilde T_1}$ and 
$\widetilde{T}_1|_{\lW_0+\ker T_1}$ are $\lH$-coercive, hence injective.

In particular, we have that the sum
\begin{equation*}
\lW_0\dotplus \ker T_1\dotplus \ker \widetilde T_1
\end{equation*}
is direct.
Indeed, let $u_0\in\lW_0$, $\nu\in\ker T_1$ and $\tilde\nu\in\ker 
\widetilde{T}_1$ be such that $u_0+\nu+\tilde\nu=0$.
Then 
$$
0=|T_1(u_0+\nu+\tilde\nu)| = |T_1(u_0+\tilde\nu)|\geq \mu_0 \|u_0+\tilde\nu\|\,,
$$
implying $u_0+\tilde\nu=0$. Acting by $\widetilde T_1$ we get
$$
0=|\widetilde{T}_1(u_0+\tilde\nu)|=|\widetilde{T}_1(u_0)|\geq \mu_0\|u_0\| \,.
$$
Thus, $u_0=0$, which implies $\tilde\nu=0$, and then finally $\nu=0$.

In the following section we shall see that
operators $T_1|_{\lW_0+\ker\widetilde T_1}$ and
$\widetilde{T}_1|_{\lW_0+\ker T_1}$ are mutually adjoint as well. 
\end{itemize}
\end{remark}

\begin{example}[Classical Friedrichs operators]\label{ex:cfo}
Let $d,r\in\N$ and $\Omega\subseteq\Rd$ be an open and bounded set with 
Lipschitz boundary $\Gamma$.
Here we present how the theory of abstract Friedrichs operators can
encompass classical Friedrichs differential operators, 
while for details we refer to \cite[Subsection 5.1]{EGC}.

We consider the restriction of operator $L$ (CFO) to $C^\infty_c(\Omega;\Cr)$ 
and denote it by $T$, i.e.
$$
T\vu =  \sum_{k=1}^d \partial_k(\mA_k \vu) + \mB \vu \,,
    \quad \vu\in C^\infty_c(\Omega;\Cr)
$$
(here the derivatives can be understood in the classical sense 
as derivatives of smooth functions are equal to their distributional
derivatives).
Since $\mB \in L^\infty(\Omega; \Mat{r})$ and 
$\mA_k \in W^{1,\infty}(\Omega;\Mat{r})$ (for any $k$),
it is obvious that $T:C^\infty_c(\Omega;\Cr)\to L^2(\Omega;\Cr)$.

For the second operator we take 
$\widetilde T:C^\infty_c(\Omega;\Cr)\to L^2(\Omega;\Cr)$ given by
$$
\widetilde T\vu =  -\sum_{k=1}^d \partial_k(\mA_k \vu) + 
    \Bigl(\mB^* + \sum_{k=1}^d \partial_k\mA_k \Bigr)\vu \,,
    \quad \vu\in C^\infty_c(\Omega;\Cr) \,.
$$
Then one can easily see that $(T,\widetilde T)$ is a joint pair of
abstract Friedrichs operators, where $\lH=L^2(\Omega;\Cr)$
and $\lD=C^\infty_c(\Omega;\Cr)$.
Indeed, (T1) is obtained by integration by parts and using (F1), 
the boundedness of coefficients implies (T2), while (T3) follows from (F2)
(a more general case where $\lH$ is taken 
to be a closed subspace of $L^2(\Omega;\Cr)$ can be found in 
\cite[Example 2]{ABCE}).

The domain of adjoint operators $T_1=\widetilde T^*$
and $\widetilde T = T^*$ (the graph space) reads
\begin{equation*}
\begin{aligned}
\lW &= \Bigl\{\vu\in L^2(\Omega;\Cr) :  
    \sum_{k=1}^d \partial_k(\mA_k \vu) + \mB \vu\in L^2(\Omega;\Cr)\Bigr\} \\
&= \Bigl\{\vu\in L^2(\Omega;\Cr) :  
    \sum_{k=1}^d \partial_k(\mA_k \vu) \in L^2(\Omega;\Cr)\Bigr\} \,.
\end{aligned}
\end{equation*}
The action of $T_1$ and $\widetilde T_1$ is (formally) the same as the action
of $T$ and $\widetilde T$, respectively (we have just that the classical 
derivatives are replaced by the distributional ones). 
It is known that $C^\infty_c(\Rd;\Cr)$ is dense in $\lW$ 
\cite[Theorem 4]{ABmc} (cf.~\cite[Chapter 1]{MJensen})
and that the boundary operator, for $\vu,\vv\in C^\infty_c(\Rd;\Cr)$,
is given by
\begin{equation*}
\dup{\lW'}{D\vu}{\vv}\lW = \int_\Gamma \mA_\mnu (\mx)\vu|_\Gamma(\mx)
    \cdot\vv|_\Gamma(\mx)\,dS(\mx) \;,
\end{equation*}
where $\mA_\mnu:=\sum_{k=1}^d\nu_k\mA_k$ and 
$\mnu=(\nu_1,\nu_2,\dots,\nu_d)\in L^\infty(\Gamma;\Rd)$ is the unit
outward normal on $\Gamma$.
In the one-dimensional case ($d=1$) for $\Omega=(a,b)$, $a<b$, 
the above formula simplifies to
\begin{equation}\label{eq:cfo-D-d1}
\dup{\lW'}{D\vu}{\vv}\lW = \mA(b)\vu(b)\cdot\vv(b)-
    \mA(a)\vu(a)\cdot\vv(a)  \,.
\end{equation}

By the definition, we have that the domain of closures
$T_0=\overline{T}$ and $\widetilde T_0=\overline{\widetilde{T}}$
is given by $\lW_0=\cl_\lW C^\infty_c(\Omega;\Cr)$, 
while by Lemma \ref{lm:boundary_op}(ii) and the identity above we have
\begin{equation*}
\begin{aligned}
\lW_0\cap C^\infty_c(\Rd;\Cr) = \Bigl\{\vu\in C^\infty_c(\Rd;\Cr) :  
    (\forall \vv &\in C^\infty_c(\Rd;\Cr)) \\
&\int_\Gamma \mA_\mnu (\mx)\vu|_\Gamma(\mx)
    \cdot\vv|_\Gamma(\mx)\,dS(\mx) = 0 \Bigr\} \;.
\end{aligned}
\end{equation*}
A more specific characterisation involving the trace operator
on the graph space can be found in \cite{ABmc, MJensen}.

This differential operator 
in the one-dimensional ($d=1$) scalar ($r=1$) 
case is the main topic of sections 
\ref{sec:1d-intro} and \ref{sec:1d-class}.
\end{example}

\begin{remark}\label{rem:intro_final}
\begin{itemize}
\item[i)] Our main concern is in applying the theory of 
abstract Friedrichs operators to classical Friedrichs 
differential operators (CFO).
However, there are many other examples that can be cast 
within this framework. 
For example in \cite[Subsection III.D]{ABCE} and 
\cite{EM19} the classical Friedrichs differential operators that are modelled 
over a closed subset of $L^2$ are presented.

Let us mention also a trivial case when $T=\widetilde T$.
For any bounded self-adjoint operator $T$ with strictly 
positive bottom, $(T,T)$ is a pair of abstract Friedrichs operators. 
Of course, the main strength of this theory is in studying 
non-symmetric operators. 

\item[ii)] Like in the classical theory, in \cite{EGC}
three different settings for imposing boundary conditions 
were introduced.
The proof that all three are equivalent was 
closed in \cite{ABcpde}.
Thus, without any loss of generality 
in this paper we chose only to work with one of them;
in terms of subspaces $\lV$ and $\widetilde \lV$. 

\item[iii)] In \cite[Lemma 8]{ABcpde} it was recognised that 
the quotient space $\widehat\lW = \lW/\lW_0$
equipped with the indefinite inner product
$$
\iscp{\hat{u}}{\hat{v}}^\wedge := \iscp uv \;,
$$
where $\hat u\equiv u+\lW_0$ is the element of $\widehat \lW$
with representative $u\in\lW$, 
is a \emph{Kre\u\i n space} (cf.~\cite{Bo}). 
This played an important role in proving the equivalence of
different abstract settings for imposing boundary conditions \cite{ABcpde}
and in obtaining the result of Theorem \ref{tm:existence-ref-op}
\cite[Theorem 13]{AEM-2017}.

\item[iv)] In previous papers on abstract Friedrichs operators, 
the characterisation of mutual adjointness from Lemma \ref{lm:mutual-adj}
was labelled as condition (V2), and this condition was set as an 
assumption e.g.~in the statements of Theorem \ref{tm:well-posedness}
and Lemma \ref{lm:coercivity}.
Here we decided to emphasise the geometric structure and 
explicitly work with mutual adjoint realisations, which is equivalent
by Lemma \ref{lm:mutual-adj}.

\item[v)] It is known that assumptions given in 
Theorem \ref{tm:well-posedness} provide only a sufficient 
condition for mutually adjoint pair 
$(T_1|_\lV, \widetilde T_1|_{\widetilde \lV})$
to be bijective (cf.~\cite[Section 6]{AEM-2017}). 
In this paper we shall focus on determining \emph{all}\/ 
bijective realisations, even if they do not come from 
that theorem, i.e.~the corresponding domains
do not satisfy condition (V1). 
\end{itemize}
\end{remark}

\section{Decomposition of the graph space}\label{sec:decomposition}

In this section we shall see that the graph space 
$\lW$ of any joint pair $(T,\widetilde T)$ of 
abstract Friedrichs operators admits a decomposition
to a direct sum (see Remark \ref{ex:V1-kernels}(ii)) 
of $\lW_0$ and kernels of adjoint operators
$T_1$ and $\widetilde T_1$. More precisely, 
the following theorem holds. 

\begin{theorem}\label{tm:decomposition}
Let $(T_0,\widetilde{T}_0)$ be a joint pair of closed
abstract Friedrichs operators on $\lH$. 
Then the following decomposition holds:
\begin{equation}\label{eq:decomposition}
    \lW \;=\; \lW_0 \dotplus \ker T_1 \dotplus \ker\widetilde T_1 \;.
\end{equation}
\end{theorem}

From the decomposition above of the graph space we can see that the \emph{room} for choosing different
boundary conditions for abstract problem $T_1 u=f$ is given by 
$\ker T_1 \dotplus \ker\widetilde T_1$.
Thus, the knowledge of adjoint operators completely describes the problem. 

Another immediate consequence of this decomposition is 
that $(T_1|_{\lW_0+\ker\widetilde T_1}, \widetilde{T}_1|_{\lW_0+\ker T_1})$
is a pair of mutually adjoint 
bijective realisations relative to $(T_0,\widetilde T_0)$.
Indeed, according to Theorem \ref{tm:well-posedness} and 
Remark \ref{ex:V1-kernels}(ii) it is sufficient to prove that 
$T_1|_{\lW_0+\ker\widetilde T_1}$ and 
$\widetilde{T}_1|_{\lW_0+\ker T_1}$ are mutually adjoint. 
For this we use Lemma \ref{lm:mutual-adj}.

\begin{corollary}\label{cor:decomposition}
Let $(T_0,\widetilde{T}_0)$ be a joint pair of closed
abstract Friedrichs operators on $\lH$. 
Then $(T_1|_{\lW_0+\ker\widetilde T_1}, \widetilde{T}_1|_{\lW_0+\ker T_1})$
is a pair of mutually adjoint 
bijective realisations relative to $(T,\widetilde T)$.
\end{corollary}

\begin{proof}
By Theorem \ref{tm:well-posedness} and Remark \ref{ex:V1-kernels}(ii)
it is sufficient to prove that
$$
\lW_0+\ker T_1 = (\lW_0+\ker\widetilde T_1)^{[\perp]}
    \qquad \hbox{and} \qquad
    \lW_0+\ker \widetilde T_1 = (\lW_0+\ker T_1)^{[\perp]} \;.
$$
Let us prove the first equality, as the proof of the second one is 
completely analogous. 

Let $u_0, v_0\in\lW_0$, $\nu\in\ker T_1$ and $\tilde{\nu}\in\ker
\widetilde T_1$ be arbitrary. Using \eqref{eq:D} and
Lemma \ref{lm:boundary_op}(i),(ii) we have
$$
\iscp{v_0+\tilde\nu}{u_0+\nu} = \overline{\iscp{\nu}{\tilde\nu}} 
    = \overline{\scp{T_1\nu}{\tilde{\nu}}} - 
    \overline{\scp{\nu}{\widetilde T_1\tilde{\nu}}} = 0 \,.
$$
Thus, $\lW_0+\ker T_1 \subseteq (\lW_0+\ker\widetilde T_1)^{[\perp]}$.

For the converse, 
let us take an arbitrary $u\in (\lW_0+\ker\widetilde T_1)^{[\perp]}$.
By \eqref{eq:decomposition} there exist 
$u_0\in\lW_0$, $\nu\in\ker T_1$ and $\tilde\nu\in\ker\widetilde T_1$
such that $u=u_0+\nu+\tilde{\nu}$.
For any $v_0\in\lW_0$ and $\tilde\nu_1\in\ker\widetilde T_1$ we have
$$
\begin{aligned}
0 = \iscp{v_0+\tilde\nu_1}{u} 
    &= \iscp{v_0+\tilde{\nu}_1}{u_0+\nu+\tilde{\nu}} \\
&= \iscp{\tilde\nu_1}{\nu} + \iscp{\tilde \nu_1}{\tilde\nu} \\
&= \iscp{\tilde \nu_1}{\tilde\nu} \,,
\end{aligned}
$$
where we have again used Lemma \ref{lm:boundary_op}(i),(ii)
together with the fact that $\ker T_1\subseteq (\ker\widetilde T_1)^{[\perp]}$. 
Putting $\tilde{\nu}_1=\tilde\nu$ we get
$$
0 = \iscp{\tilde \nu}{\tilde\nu} = \scp{T_1\tilde{\nu}}{\tilde\nu}
= \scp{(T_1+\widetilde T_1)\tilde{\nu}}{\tilde\nu} 
    \geq 2\mu_0 \|\tilde\nu\|^2 \,,
$$
where the last inequality is due to condition (T3). 
Hence, necessarily $\tilde{\nu}=0$, which implies 
$u=u_0+\nu\subseteq \lW_0+\ker T_1$.
\end{proof}

\begin{remark}
\begin{itemize}
\item[i)] It is easy to see that the property that 
$(T_1|_{\lW_0+\ker\widetilde T_1}, \widetilde{T}_1|_{\lW_0+\ker T_1})$
is a pair of mutually adjoint 
bijective realisations relative to $(T,\widetilde T)$
is in fact equivalent to the decomposition \eqref{eq:decomposition}.
Indeed, we just consider $T_1|_{\lW_0+\ker\widetilde T_1}$
as the reference operator and apply the decomposition of the graph space
given in point (i) of Appendix to retrieve \eqref{eq:decomposition}.

Using this equivalence one can construct an alternative proof of 
the decomposition \eqref{eq:decomposition} from the one presented in the rest
of the section.

\item[ii)] Another equivalent statement of Theorem 
\ref{tm:decomposition} can be formulated in terms 
of the Kre\v\i n space $\widehat\lW$
(see Remark \ref{rem:intro_final}(iii)).
Indeed, \eqref{eq:decomposition} holds
if and only if subspaces $\widehat{\ker T_1}$
and $\widehat{\ker\widetilde T_1}$ provide a canonical
(or fundamental) decomposition of $\widehat \lW$
(cf.~\cite{Bo}).
\end{itemize}
\end{remark}

\begin{remark}
The previous corollary implies that for any 
pair of abstract Friedrichs operators $(T,\widetilde T)$
there exists a pair $(T_\mathrm{r},
T_\mathrm{r}^*)$ of mutually adjoint 
bijective realisations such that 
$\dom T_\mathrm{r}+\dom T_\mathrm{r}^*$ is 
closed in $\lW$ and
$\dom T_\mathrm{r}+\dom T_\mathrm{r}^*=\lW$.
This situation is beneficial for studying 
other descriptions of boundary conditions
(see \cite[Section 4]{EGC} and \cite[Section 3]{ABcpde}).
\end{remark}

Let us go back to the decomposition \eqref{eq:decomposition}. 
In order to prove it, we shall first prove several 
auxiliary claims. 
In all of them, as well as in the rest of the section, 
we assume that we are given a joint pair 
$(T_0,\widetilde T_0)$ of closed abstract Friedrichs operators on $\lH$.
We start by proving that 
the sum $\lW_0 \dotplus \ker{T_1} \dotplus \ker{\widetilde T_1}$
is closed in $\lW$. 

\begin{lemma}\label{lm:decomp-closed}
The sum $\lW_0 \dotplus \ker{T_1} \dotplus \ker{\widetilde T_1}$
is direct and closed in $\lW$.

In particular, $\lW_0 \dotplus \ker{T_1}$ and $\lW_0 \dotplus \ker{T_1}$ 
are both closed in $\lW$.
\end{lemma}
\begin{proof}
The second part of the statement is just a simple consequence, 
so let us just focus on studying 
$\lW_0 + \ker{T_1} + \ker{\widetilde T_1}$.

The fact that the sum is direct has already been shown in Remark 
\ref{ex:V1-kernels}(ii).

Let us take an arbitrary convergent sequence $u_n=u_n^0+\nu_n+\tilde\nu_n$ 
in $\lW_0 \dotplus \ker{T_1} \dotplus \ker{\widetilde T_1}$
($u_n^0\in\lW_0$, $\nu_n\in\ker T_1$, $\tilde\nu_n\in\ker \widetilde{T}_1$)
with respect to the graph norm
and denote by $u\in\lW$ its limit.

Since $T_1(u_n^0+\nu_n+\tilde\nu_n) = T_1(u_n^0+\tilde\nu_n)$
is a Cauchy sequence in $\lH$ and 
$T_1|_{\lW_0+\ker\widetilde{T}_1}$ is $\lH$-coercive (see Remark
\ref{ex:V1-kernels}(ii)), $(u_n^0+\tilde\nu_n)$ is a Cauchy sequence in 
$\lH$ as well, hence convergent. 
Let us define $w:=\lim_n (u_n^0+\tilde\nu_n)\in\lH$ and 
$\nu:=u-w\in\lH$.

From
\begin{align*}
    \|\nu_n-\nu\| & =\|(u_n^0+\nu_n+\tilde\nu_n)-u-(u_n^0+\tilde\nu_n-w)\| \\
    & \leq \|u_n^0+\nu_n+\tilde\nu_n-u\|+\|u_n^0+\tilde\nu_n-w\| \,,
\end{align*}
we have $\lim_n \nu_n = \nu$. Since $\ker{T_1}$ is closed in $\mathcal{H}$, this implies that $\nu\in\ker T_1$ and
$\nu_n\xrightarrow{\mathcal{W}}\nu$.

Therefore, we have so far that $u_n^0+\tilde\nu_n\xrightarrow{\mathcal{W}} u-\nu$, implying that 
$\widetilde T_1(u_n^0+\tilde\nu_n)=\widetilde T_0(u_n^0)$
is a Cauchy sequence in $\lH$. 
Since $\widetilde T_0$ is $\lH$-coercive (see Remark
\ref{ex:V1-kernels}(i)), $(u_n^0)$ is also a Cauchy sequence in 
$\lH$, hence convergent in $\lH$.
By the closedness of $\widetilde T_0$ this implies that in fact 
$(u_n^0)$ is convergent in $\lW_0$ (in the graph norm) and let us denote 
its limit by $u_0\in\lW_0$.

Let us define $\tilde\nu:=u-u_0-\nu$.
Analogously as for $(\nu_n)$, we get that $\tilde\nu_n\xrightarrow{\mathcal{W}}\tilde\nu\in \ker\widetilde T_1$.
Thus, $u_n^0+\nu_n+\tilde\nu_n\xrightarrow{\mathcal{W}} u_0+\nu+\tilde\nu$.
Uniqueness of the limit finally implies 
$u=u_0+\nu+\tilde\nu\in 
\lW_0 \dotplus \ker{T_1} \dotplus \ker{\widetilde T_1}$.
\end{proof}

The previous lemma together with Lemma \ref{lm:closed-compl}
implies that
\begin{equation*}
\lW_0 \dotplus \ker T_1 \dotplus \ker \widetilde T_1 \;=\;
    \bigl(\lW_0 \dotplus \ker T_1 \dotplus 
    \ker \widetilde T_1\bigr)^{[\perp][\perp]} \;.
\end{equation*}
Thus, in order to prove Theorem \ref{tm:decomposition}
it is sufficient to show the following equality:
$$
\bigl(\lW_0 \dotplus \ker T_1 \dotplus 
    \ker \widetilde T_1\bigr)^{[\perp][\perp]}\;=\;\lW \;.
$$
That will be obtained using the following refinement of the 
decomposition given in part (i) of Appendix.

\begin{lemma}\label{lm:grubb-decomp}
For any bijecitve realisation $T_\mathrm{r}$ of $T_0$ 
(see p.~\pageref{eq:realisations}), we have
\begin{align*}
\mathcal{W}  & \;=\; \mathcal{W}_0 \,\dotplus\, 
    T_r^{-1}(\ker \widetilde{T}_1) \,\dotplus\, \ker T_1 \\
& \;=\; \mathcal{W}_0 \,\dotplus\, ({T_r^*})^{-1}(\ker{T}_1)
    \,\dotplus\, \ker \widetilde{T}_1 \;.
\end{align*}
\end{lemma}

\begin{proof}
It is evident that two given decompositions are symmetrical
(we get the second one by replacing the role of $T_0$ and $\widetilde T_0$),
so we present a proof only for the first one.

By the decomposition given in part (i) of Appendix we have
\begin{equation*}
\lW \;=\; \dom T_\mathrm{r} \dotplus \ker T_1 \;.
\end{equation*}
Thus, it is sufficient to prove
\begin{equation}\label{eq:domTr-decomp}
\dom T_\mathrm{r}\;=\;\lW_0\dotplus T_\mathrm{r}^{-1}(\ker \widetilde{T}_1)
\end{equation}
(we refer to page \pageref{eq:realisations} where existence and 
properties of $T_\mathrm{r}^{-1}$ were discussed). 

Since $T_0\subseteq T_\mathrm{r}$, we have $\lW_0\subseteq \dom T_\mathrm{r}$,
while inclusion $T_\mathrm{r}^{-1}(\ker \widetilde{T}_1)\subseteq 
\dom T_\mathrm{r}$ is trivial. 
Hence, 
$\lW_0 + T_\mathrm{r}^{-1}(\ker \widetilde{T}_1)\subseteq 
\dom T_\mathrm{r}$.
Now let $u$ be an arbitrary element in $\dom T_\mathrm{r}$.
Since $T_\mathrm{r} u\in\lH$, by \eqref{eq:orth-decomp-lH}
there exist $u_0\in\lW_0$ and $\tilde\nu\in\ker \widetilde T_1$
such that $T_\mathrm{r} u = T_0 u_0 + \tilde\nu
=T_\mathrm{r}u_0+\tilde\nu$.
Thus, using that $T_1|_{\dom T_\mathrm{r}}=T_\mathrm{r}$ is a bijection, 
we have
\begin{align*}
u =T_\mathrm{r}^{-1}T_\mathrm{r}u &= T_r^{-1}(T_\mathrm{r}u_0+\tilde\nu) \\
     & = u_0 + T_\mathrm{r}^{-1}(\tilde\nu) \;,
\end{align*}
implying $u\in\lW_0 + T_\mathrm{r}^{-1}(\ker \widetilde{T}_1)$.

It is left to see that the sum $\lW_0 + T_\mathrm{r}^{-1}
(\ker \widetilde{T}_1)$ is direct. 
Let $u_0\in\lW_0$ and $\tilde\nu\in\ker \widetilde T_1$
such that $u_0+T_\mathrm{r}^{-1}(\tilde\nu)=0$.
Acting with $T_\mathrm{r}$ we get
$\tilde\nu = -T_0u_0$. Hence, $\tilde\nu\in\ker\widetilde T_1\cap
\ran T_0 =\{0\}$, implying $\tilde\nu=T_0u_0=0$.
Since $T_0$ is injective, we have $u_0=0$ as well.
This establishes equality \eqref{eq:domTr-decomp}
with which the proof is completed.
\end{proof}

Let us now go back to the above mentioned sufficient condition. 

\begin{lemma}\label{lm:orth-orth}
It holds $\bigl(\lW_0 \dotplus \ker T_1 \dotplus 
    \ker \widetilde T_1\bigr)^{[\perp][\perp]}\;=\;\lW \;.$
\end{lemma}
\begin{proof}
Since $\lW_0^{[\perp]} = \lW$, it is sufficient to prove
$$
\bigl(\lW_0 \dotplus \ker T_1 \dotplus 
    \ker \widetilde T_1\bigr)^{[\perp]}\;=\;\lW_0 \;.
$$
Since $\ker D=\lW_0$, inclusion 
$\lW_0\subseteq\bigl(\lW_0 \dotplus \ker T_1 \dotplus 
    \ker \widetilde T_1\bigr)^{[\perp]}$ is obvious.
    
Let us take an arbitrary $u\in \bigl(\lW_0 \dotplus \ker T_1 \dotplus 
    \ker \widetilde T_1\bigr)^{[\perp]}$
and let $T_\mathrm{r}$ be a bijective realisation of $T_0$
(it exists by Theorem \ref{tm:existence-ref-op}(i)).
Since $u\in\lW$, by Lemma \ref{lm:grubb-decomp} there exist
unique $u_0\in\lW_0$, $\nu\in\ker T_1$ and $\tilde\nu\in\ker
\widetilde T_1$ such that
$u=u_0+T_\mathrm{r}^{-1}(\tilde\nu)+\nu$.
For arbitrary $v_0\in\lW_0$, $\nu_1\in\ker T_1$ and $\tilde\nu_1\in\ker
\widetilde T_1$ we have
\begin{equation}\label{eq:temp_decomp}
\begin{aligned}
0=\iscp{u}{v_0+\nu_1+\tilde\nu_1} 
    & = \iscp{u_0+T_\mathrm{r}^{-1}(\tilde\nu)+\nu}{v_0+\nu_1+\tilde\nu_1} \\
& = \iscp{T_\mathrm{r}^{-1}(\tilde\nu)+\nu}{\nu_1+\tilde\nu_1} \\
& = \iscp{T_\mathrm{r}^{-1}(\tilde\nu)}{\nu_1} 
    + \iscp{T_\mathrm{r}^{-1}(\tilde\nu)}{\tilde\nu_1}
    + \iscp{\nu}{\nu_1}
    + \iscp{\nu}{\tilde\nu_1} \\
& = \iscp{T_\mathrm{r}^{-1}(\tilde\nu)}{\nu_1} 
    + \iscp{T_\mathrm{r}^{-1}(\tilde\nu)}{\tilde\nu_1}
    + \iscp{\nu}{\nu_1} \;,
\end{aligned}
\end{equation}
where in the third equality we have used that 
$\lW_0^{[\perp]}=\lW$ (see Lemma \ref{lm:boundary_op}) 
and in the last that 
$\ker T_1\subseteq \bigl(\ker \widetilde T_1\bigr)^{[\perp]}$.
Now for $\nu_1=0$ and $\tilde\nu_1=\tilde\nu$ we get
$$
0 = \iscp{T_\mathrm{r}^{-1}(\tilde\nu)}{\tilde\nu}
    = \scp{\tilde\nu}{\tilde\nu} = \|\tilde\nu\|^2 \,,
$$
where in the second equality we have used 
the fact that $T_\mathrm{r}\subseteq T_1$ and that 
$\tilde\nu\in\ker \widetilde T_1$ (see \eqref{eq:D}).
Thus, $\tilde\nu=0$.

Returning to \eqref{eq:temp_decomp} with $\nu_1=\nu$
and using that $\widetilde T_1|_{\lW_0+\ker T_1}$
is $\lH$-coercive (see Remark \ref{ex:V1-kernels}(ii)),
we get 
$$
0 = |\iscp\nu\nu| = |\scp{\widetilde T_1\nu}{\nu}| \geq 
    \mu_0 \|\nu\|^2 \,.
$$
Therefore, $\nu=0$ as well. 

This implies that $u=u_0\in\lW_0$,
which was to be shown.
\end{proof}

Now we are ready to prove Theorem \ref{tm:decomposition}.

\begin{proof}[Proof of Theorem \ref{tm:decomposition}]
By Lemma \ref{lm:decomp-closed} we have that the sum
$\lW_0 \dotplus \ker T_1 \dotplus \widetilde T_1$
is direct and closed in $\lW$.
Then applying Lemma \ref{lm:closed-compl} we get
that the following equality holds:
\begin{equation*}
\lW_0 \dotplus \ker T_1 \dotplus \ker \widetilde T_1 \;=\;
    \bigl(\lW_0 \dotplus \ker T_1 \dotplus 
    \ker \widetilde T_1\bigr)^{[\perp][\perp]} \;.
\end{equation*}
Finally, using Lemma \ref{lm:orth-orth} we reach to 
\eqref{eq:decomposition}.
\end{proof}

\begin{remark}
In the case of finite dimensional kernels, 
i.e.~$\dim \ker T_1<\infty$ and $\dim \ker\widetilde T_1<\infty$,
the statement of Theorem \ref{tm:decomposition} is a direct consequence 
of Lemma \ref{lm:grubb-decomp}.
Indeed, by Lemma \ref{lm:grubb-decomp} we get
$\lW/\lW_0  \cong T_\mathrm{r}^{-1}(\ker \widetilde{T}_1)\dotplus
\ker T_1$.
Since $T_r :\dom T_r\to \mathcal{H}$ is a bijection, we have
\begin{align*}
\dim(\lW/\lW_0) &= \dim\bigl(T_\mathrm{r}^{-1}(\ker 
    \widetilde{T}_1)\dotplus\ker T_1\bigr) \\
&= \dim(\ker\widetilde T_1) +\dim(\ker T_1) <\infty \,,
\end{align*}
hence the codimension is finite.

On the other hand, obviously $(\lW_0\dotplus\ker T_1\dotplus
\ker\widetilde T_1)/\lW_0\subseteq \lW/\lW_0$
and (since the sum is direct)
\begin{align*}
\dim\bigl((\lW_0\dotplus\ker T_1\dotplus \ker\widetilde T_1)/\lW_0\bigr)
    &= \dim(\ker T_1) + \dim(\ker \widetilde T_1) \\
    &= \dim(\lW/\lW_0) \;.
\end{align*}
Therefore, we get $(\lW_0\dotplus\ker T_1\dotplus \ker\widetilde T_1)/\lW_0 = \lW/\lW_0$, implying \eqref{eq:decomposition}.
\end{remark}

\begin{remark}
    The decomposition \eqref{eq:decomposition} can be seen as a von Neumann type formula for non-symmetric operators (see 
    e.g.~\cite[Proposition 3.7]{Schmudgen}). At first glance, the required assumption of coercivity ((T3) condition) might look unsatisfactory, as it is not assumed in the symmetric case. However, when looking into a standard proof of the von Neumann formula, one can see that it is also used there, but it was guaranteed by the fact that for a symmetric densely defined operator $S$ and $\lambda\in \C\setminus \R$ we have that $S-\lambda \mathbbm{1}$ is coercive (see \cite[Proposition 3.2(i)]{Schmudgen}). Moreover, this property is then used in principle for the same reasons as here (e.g.~closedness of $\ran T_0$, injectivity of certain operators).
    
    Let us emphasise that the von Neumann formula cannot be derived from 
    Theorem \ref{tm:decomposition}: if for a symmetric operator $S$ the pair $(S-\lambda\mathbbm{1}, S-\bar{\lambda}\mathbbm{1})$ is a joint pair of abstract Friedrichs operators, $S$ is necessarily bounded. 
    Thus, it might be of an independent interest to study minimal requirements
    on dual pairs of operators \eqref{eq:dual-pairs} 
    for which the decomposition \eqref{eq:decomposition} holds.
\end{remark}

For any bijective realisation of $T_0$, the  decomposition given in part (i) of Appendix holds (which was used in the proof of Lemma  \ref{lm:grubb-decomp}). In fact we shall see in the following lemma that the opposite implication holds as well. 

\begin{lemma}\label{lm:con}
 Let $\lV $ be a closed subspace  of  $\lW$ such that $\lW_0\subseteq \lV$. Then $T_1|_{\lV}$ is bijective
 if and only if\/ $\lV\dot{+}\ker T_1=\lW$.
\end{lemma}
\begin{proof}
The first implication is followed from the decomposition given in part (i) of Appendix.

For the converse,
$\lV\cap \ker T_1=\{0\}$ implies that $T_1|_{\lV}$ is injective. Now let $f\in \lH$. Since $T_1:\lW\to \lH$ is surjective, there exists $u\in \lW$ such that $T_1u=f$. We also have for some $v\in \lV$ and $\nu \in \ker T_1$, $u=v+\nu$. Thus,
\begin{align*}
    f=T_1u=T_1(v+\nu)=T_1v=T_1|_{\lV}
\end{align*}
implies that $T_1|_{\lV}$ is surjective. Hence, $T_1|_{\lV}$ is a bijective realisation.
\end{proof}

\begin{remark}
\begin{itemize}
    \item[i)] Since $\lW_0\dot{+}\ker \widetilde{T}_1$ is a closed subspace of $\lW$ with $\lW_0\subseteq \lW_0\dot{+}\ker \widetilde{T}_1$, we have that $T_1|_{\lW_0\dot{+}\ker \widetilde{T}_1}$ is a bijective realisation of $T_0$. 
    This allows us to construct an alternative proof of Corollary 3.2.
    
    \item[ii)] For any closed subspace  $\lV$ of  $\lW$ such that $\lW_0\subseteq \lV$ and $T_1|_{\lV}$ is bijective  we have $\lV/\lW_0\cong \, \ker \widetilde{T}_1$.
\end{itemize}

\end{remark}

\section{One-dimensional scalar (CFO):  
    Preliminaries}\label{sec:1d-intro}

In this section and in the rest of the manuscript we study 
(CFO) in the one-dimensional ($d=1$) scalar ($r=1$) case. 
For the domain we take an open interval $\Omega=(a,b)$,
$a<b$. Then $\lD=C^\infty_c(a,b)$ and $\lH=L^2(a,b)$.
We adjust the notation of $T, \widetilde T:\lD\to\lH$
given in Example \ref{ex:cfo} in the following way:
\begin{equation}\label{eq:TTtilda}
T\phi := (\alpha \phi )'+\beta \phi
	\qquad \hbox{and} \qquad
	\widetilde{T}\phi:=-(\alpha \phi)'
	+(\overline{\beta}+\alpha ')\phi \;,
\end{equation}
where $\alpha\in W^{1,\infty}((a,b);\R)$,
$\beta\in L^\infty((a,b);\C)$ and for some $\mu_0>0$
we have $2\Re \beta +\alpha '\geq 2\mu_0>0$
($\Re z$ denotes the real part of complex number $z$
and $'$ the derivative). 

It is commented in Example \ref{ex:cfo} that 
$(T, \widetilde T)$ 
is a joint pair of abstract Friedrichs operators.
Moreover, the graph space is given by 
\begin{equation*}
\lW \;=\; \bigl\{ u\in\lH : (\alpha u)'\in\lH\bigr\} \;,
\end{equation*}
while the graph norm is equivalent to 
$\|u\|_{\mathcal{W}}:=\|u\|+\|(\alpha u)'\|$
($\|\cdot\|$ stands, as usual, for the norm on $\lH$
induced by the standard inner product, i.e.~the $L^2$ norm
on $(a,b)$).
In fact, $u\in\lH$ belongs to $\lW$ if and only if 
$\alpha u\in H^1(a,b)$. Thus, by the standard Sobolev 
embedding theorem (see e.g.~\cite[Theorem 8.2]{Brezis})
for any $u\in\lW$ we have $\alpha u\in C([a,b])$.
This in particular implies that for any $u\in\lW$ and
$x\in[a,b]$
evaluation $(\alpha u)(x)$ is well defined.
However, $\alpha(x)u(x)$ is not necessarily
meaningful as $u$ itself is not necessarily continuous.
A more precise description of the graph space 
is given in the following lemma.

\begin{lemma}
    Let $I:=[a,b]\setminus\alpha^{-1}(\{0\})$.
    Then $\lW\subseteq H^1_\mathrm{loc}(I)$,
    i.e.~for any $u\in\lW$ and $[c,d]\subseteq I$, $c<d$,
    we have $u|_{[c,d]}\in H^1(c,d)$.
\end{lemma} 
\begin{proof} 
Since $\alpha$ is continuous, $I$ is relatively open in $[a,b]$.
Let us take $[c,d]\subseteq I$, $c<d$ (if such segment does not exist, 
then $\alpha\equiv 0$ and $I=\emptyset$, which is a trivial case), 
and define
$\alpha_0:=\min_{x\in[c,d]} |\alpha(x)|$.
Obviously $\alpha_0>0$. 

Let $u\in {C}_c^{\infty}(\mathbb{R})$, 
then
\begin{align*}
\|u'\|_{{L}^2(c,d)} 
&\leq\frac{1}{\alpha_0}\|\alpha u'\|_{{L}^2(c,d)}\\
&\leq\frac{1}{\alpha_0}\Bigl(\|(\alpha u)'\|_{{L}^2(c,d)}
    +\|\alpha'u\|_{{L}^2(c,d)}\Bigr)\\
&\leq\frac{1}{\alpha_0}\Bigl(\|(\alpha u)'\|
    +\|\alpha'\|_{{L}^{\infty}(a,b)}\|u\| \Bigr)
    \leq \frac{1+\|\alpha\|_{W^{1,\infty}(a,b)}}{\alpha_0}
    \|u\|_\lW \;.
\end{align*}
Since $C^\infty_c(\R)$ is dense in $\lW$ (cf.~\cite[Theorem 4]{ABmc}), 
by a standard argument we can deduce that $u|_{[c,d]}\in H^1(c,d)$
and there exists $C>0$ (dependent on $c, d$) such that
\begin{equation*}
\|u\|_{H^1(c,d)}\leq C \|u\|_\lW \;, \quad u\in\lW \;.
\end{equation*}
\end{proof}

\begin{remark}\label{rem:h1loc}
\begin{itemize}
    \item[i)] If $x\in I$, where $I$ is defined in the statement of the
    previous lemma, then we can write $(\alpha u)(x)=\alpha(x) u(x)$.
    Then it is natural to extend $(\alpha u)(x)$ to be 0 if $x\not\in I$.
    \item[ii)] From the proof of the previous lemma we can deduce 
    that $\lW$ is continuously embedded in $H^1_\mathrm{loc}(I)$,
    which we write as $\lW\hookrightarrow H^1_\mathrm{loc}(I)$.
    Of course the embedding is strict, which we illustrate on examples
    in Section \ref{sec:examples}.
    This can also be argued 
    by noting that $H^1_\mathrm{loc}(I)$ is not a normed space
    (it is a Fr\'echet space), while $\lW$ is.
    
    On the other hand, it is trivial to see that 
    $H^1(a,b)\hookrightarrow \lW$ and $H^1(a,b)=\lW$ if and only if
    $\alpha$ has no zeros on $[a,b]$.
\end{itemize}
\end{remark}

In \eqref{eq:cfo-D-d1} an explicit formula for the boundary operator 
$D$ is given on the dense subspace $C^\infty_c(\R)$.
By the previous lemma we can extend it uniquely
on $\lW$ by density argument, which reads
\begin{equation}\label{eq:D-scalar1d}
\dup{\lW'}{Du}{v}\lW = \bigl(\alpha u\overline v\bigr)(b) -
    \bigl(\alpha u\overline v\bigr)(a) \;, \quad u,v\in\lW\;,
\end{equation}
where we define (see Remark \ref{rem:h1loc}(i))
\begin{equation}\label{eq:alphau}
\bigl(\alpha u\overline{v}\bigr)(x) :=
    \left\{\begin{array}{lcl}0&,& \alpha(x)=0\\
        \alpha(x)u(x)\overline{v(x)} &,& \alpha(x)\neq 0 
        \end{array}\right. \;, \quad x\in [a,b]\,.
\end{equation}

The domain of the closures $T_0$ and $\widetilde T_0$
satisfies $\lW_0=\cl_\lW C^\infty_c(\R)$,
but having \eqref{eq:D-scalar1d} it is easier to use 
$\ker D=\lW_0$ to characterise $\lW_0$.

\begin{lemma}\label{lm:W0}
The space $\mathcal{W}_0$ can be characterised as
\begin{equation*}
\lW_0=\Bigl\{u\in \lW : 
    (\alpha u)(a)=(\alpha u)(b) =0 \Bigl\} \;,
\end{equation*}
where $(\alpha u)(x)$ is to be understood as in
\eqref{eq:alphau} (see also Remark \ref{rem:h1loc}(i)).
\end{lemma}
\begin{proof} 
We know that $\lW_0=\ker D$ (see Lemma \ref{lm:boundary_op}(ii)). Let $u \in \lW$ be such that 
$(\alpha u)(a)=(\alpha u)(b)=0$.
Then for any $v \in \lW$ by \eqref{eq:D-scalar1d} we have
\[
\dup{\lW'}{Du}{v}\lW = 
    (\alpha u \overline{v})(b)
    -(\alpha u \overline{v})(a) =0-0=0
\]
(of course we deal with this in cases depending on the value of $\alpha$ at boundary points, but in each case the result is 0). 
Hence,
$$
\bigl\{u\in \lW : 
    (\alpha u)(a)=(\alpha u)(b) =0 \bigl\}\subseteq \ker D \,.
$$

Conversely, let $u \in \ker D\subseteq\lW$.
Then for any $v\in H^1(a,b)\subseteq\lW$ it holds
\begin{align*}
0 = \dup{\lW'}{Du}{v}\lW &=
    (\alpha u\overline{v})(b) - (\alpha u\overline{v})(a) \\
&= (\alpha u)(b)\overline{v(b)} - (\alpha u)(a)\overline{v(a)} \;, 
\end{align*}
where we have used that $v$ is continuous.
For $v(x)=x-a$ we get $(\alpha u)(b)=0$,
while for $v(x)=x-b$ we reach to $(\alpha u)(a)=0$,
completing the proof.
\end{proof}

From the decomposition \eqref{eq:decomposition}
we see that $\ker T_1+\ker \widetilde T_1$, or equivalently
$\lW/\lW_0$, plays an important role in studying boundary conditions
associated to $T$ (or $\widetilde T$).
Here we present a result on the codimension.
        
\begin{lemma}\label{lm:codim}
 \begin{equation*}
 \dim(\mathcal{W}/{\mathcal{W}_0)}\,=\;
    \left\{\begin{array}{lcl}
                       2 &,& \alpha(a)\alpha(b)\neq 0 \;, \\ 
                       1 &,& \bigl(\alpha(a)=0 \wedge 
                        \alpha(b)\neq 0\bigr)\,\vee \, 
                        \bigl(\alpha(a)\neq 0 \wedge 
                         \alpha(b)=0\bigr)\;, \\
                       0 &,& \alpha(a)=\alpha(b)=0 \;.
         \end{array}\right.
\end{equation*}
\end{lemma}         
\begin{proof} 
Choose $\phi , \psi \in H^1(a,b)$, such that 
$\phi(a)=1$, $\phi(b)=0$ and 
$\psi(a)=0$, $ \psi(b)=1$.
Define $\hat{\phi}:=\phi+\mathcal{W}_0$ and $\hat{\psi}:=\psi+\mathcal{W}_0$. 

For $\alpha(a)\alpha(b)\neq 0$ we prove that the set 
$E:=\{\hat{\phi},\hat{\psi}\}$ is a basis of $\mathcal{W}/{\mathcal{W}_0}$.
First we claim that the set $E$ is linearly independent. 
Indeed, if $E$ were linearly dependent then for some 
non-zero scalar $r$ we would have $\hat{\psi}=r\hat{\phi}$,
implying $\hat{\psi}-r\hat{\phi}=\hat{0}=\mathcal{W}_0$.
Hence, $\psi-r\phi \in \mathcal{W}_0$, which by Lemma \ref{lm:W0}
leads to
$
\bigl(\alpha (\psi-r\phi)\bigr)(a)
    =\bigl(\alpha (\psi-r\phi)\bigr)(b)=0
$.
But, 
\[
\bigl(\alpha (\psi-r\phi)\bigr)(a)
    =\alpha(a)\psi(a)-r\alpha(a)\phi(a)=-r\alpha(a)\neq 0 \,,
\]
which is a contradiction. Hence, 
the set $E$ is linearly independent.
Now let $u\in \mathcal{W}$, then 
\[
u-u(a)\phi -u(b)\psi \in \mathcal{W}_0 \;,
\] 
which means that $E$ spans $\mathcal{W}/{\mathcal{W}_0}$.
Therefore, $E$ is a basis of $\mathcal{W}/{\mathcal{W}_0}$ and $\dim(\mathcal{W}/{\mathcal{W}_0})=2$.

If $\alpha(a)=0$ and $\alpha(b)\neq 0$, 
then $\phi \in \mathcal{W}_0$, so $\mathcal{W}/\mathcal{W}_0=\operatorname{span}\{\hat{\psi}\}$ and $\dim(\mathcal{W}/\mathcal{W}_0)=1$.
Similarly, if $\alpha(a)\neq 0$ and $\alpha(b)=0$, we also 
have $\dim(\mathcal{W}/\mathcal{W}_0)=1$.

If $\alpha(a)=\alpha(b)=0$, then $D=0$, hence 
$\mathcal{W}=\ker(D)=\mathcal{W}_0$,
implying $\dim(\mathcal{W}/\mathcal{W}_0)=0$.
\end{proof}

\begin{remark}\label{rem:codim}
\begin{itemize}
    \item[i)] If $\min_{x\in[a,b]}|\alpha(x)|>\alpha_0>0$
    (see Remark \ref{rem:h1loc}(ii)), 
    then the statement of the previous lemma 
    reveals a well known fact that 
    $\dim \bigl(H^1(a,b)/H^1_0(a,b)\bigr) = 2$.
    \item[ii)] By the decomposition \eqref{eq:decomposition}
    we have
    $$
    \dim(\ker T_1) + \dim(\ker\widetilde T_1) = 
       \dim \lW/\lW_0 \,.
    $$
    Thus, by the previous lemma and 
    Theorem \ref{tm:existence-ref-op} we can 
    immediately conclude that in the case 
    $\alpha(a)\alpha(b)=0$ there is only one bijective 
    realisation of $T_0$. Moreover, in the opposite 
    case $\alpha(a)\alpha(b)\neq 0$ there are infinitely
    many bijective realisations if and only if
    $\dim(\ker T_1) = \dim (\ker\widetilde{T}_1)$.
    
    We shall justify and improve these conclusions by a direct 
    inspection in the following section.
\end{itemize}
\end{remark}

\section{One-dimensional scalar (CFO):
    Classification}\label{sec:1d-class}
 
This section contains a complete classification of
all bijective realisations relative to the pair
\eqref{eq:TTtilda} from the previous section.
More precisely, we seek for all pairs of subspaces 
$(\lV,\widetilde\lV)$ such that 
$(T_1|_\lV, \widetilde T_1|_{\widetilde\lV})$
is a pair of mutually adjoint \emph{bijective} realisations
relative to $(T,\widetilde T)$, 
where $T$ and $\widetilde T$ are given by
\eqref{eq:TTtilda}.
The analysis is divided into three cases depending on $\operatorname{sign}\bigl(\alpha(a)\alpha(b)\bigr)$.

\subsection{Case 1: $\alpha(a)\alpha(b)=0$}

Let us consider first the subcase $\alpha(a)=\alpha(b)=0$.
Then the boundary map is trivial, i.e.~$D=0$. 
This implies $\lW_0=\ker(D) = \lW$,
thus the only possible choice is $(\mathcal{V},\widetilde{\mathcal{V}})=(\mathcal{W},\mathcal{W})$.

Now we deal with the subcase in which exactly one of numbers
$\alpha(a)$, $\alpha(b)$ is equal to zero. 
Let us present in a full detail only the situation where
$\alpha(a)=0$ and $\alpha(b)>0$, as analysis 
for the other three is completely analogous.
Moreover, one can derive results for the other situations by switching 
the role of $T$ and $\widetilde T$ and/or reflecting 
operators, i.e.~changing the domain to $(-b,-a)$.

The boundary map in the case $\alpha(a)=0$ and $\alpha(b)>0$ 
reads
\[ 
\dup{\lW'}{Du}{v}\lW = \alpha(b)u(b)\overline{v(b)} \;,
    \quad u,v\in\lW \;,
\]
and (see Lemma \ref{lm:W0}) 
$\mathcal{W}_0=\{u\in \mathcal{W}:u(b)=0\}$.
Since for any $u\in\lW$ we have
\begin{equation}\label{eq:case1}
\dup{\lW'}{Du}{u}\lW = \alpha(b)|u(b)|^2 \geq 0 \;,
\end{equation}
pair $(\lW,\lW_0)$ satisfies condition (V1). 
Furthermore, $(T_1|_{\lW},\widetilde T_1|_{\lW_0})
=(T_1,\widetilde T_0)$ is trivially a pair of 
mutually adjoint operators. 
Therefore, by Theorem \ref{tm:well-posedness}
this pair forms a pair of mutually adjoint bijective 
realisations relative to $(T,\widetilde T)$.
Since this implies that $\ker T_1=\{0\}$,
by Theorem \ref{tm:existence-ref-op}(ii), 
$(T_1,\widetilde T_0)$ is the only pair 
of mutually adjoint bijective 
realisations relative to $(T,\widetilde T)$.

The overall conclusion is: 
\begin{equation}\label{eq:case1-VVtilda}
(\lV,\widetilde\lV) \;=\;
    \left\{\begin{array}{lcl}
        (\lW,\lW_0) &,& \bigl(\alpha(a)=0 \wedge \alpha(b)\geq 0)
            \,\vee\, \bigl(\alpha(a)\leq 0 \wedge \alpha(b)=0 \bigr)\\
        (\lW_0,\lW) &,& \bigl(\alpha(a)=0 \wedge \alpha(b)\leq 0)
            \,\vee\, \bigl(\alpha(a)\geq 0 \wedge \alpha(b)=0 \bigr)\end{array}\right. \;,
\end{equation}
i.e we always have only one pair of bijective realisation.
In the above we also included the first case $\alpha(a)=\alpha(b)=0$
as then $(\lW,\lW_0)=(\lW_0,\lW)=(\lW,\lW)$.

With this analysis we obtained the same conclusion for this case
as in Remark \ref{rem:codim}(ii), but without using
Lemma \ref{lm:codim}.
Although we have fully characterised bijective realisations, 
let us say a little more about kernels of $T_1$ and $\widetilde T_1$.

In the case $\alpha(a)=\alpha(b)=0$ it is clear that
$\ker T_1 =\ker \widetilde T_1 = \{0\}$.
This means that both equations
\begin{equation*}
(\alpha\phi)'+\beta\phi = 0 \quad \hbox{and} \quad
    -(\alpha\phi)' + (\overline{\beta}+\alpha')\phi=0
\end{equation*}
do not have any non-trivial solution in $\lW$.

If exactly one of numbers $\alpha(a)$ and $\alpha(b)$
is equal to zero, from Remark \ref{rem:codim}(ii) we have
$\dim(\ker T_1) + \dim(\ker \widetilde T_1) =1$,
while the analysis above implies that one of dimensions 
equals 0 (the one associated to the operator for which we took 
the whole graph space $\lW$ as the domain of the bijective 
realisation -- see \eqref{eq:case1-VVtilda}). 
To be more specific, let us stick to the case 
$\alpha(a)=0$ and $\alpha(b)>0$.
Then, $\dim(\ker T_1)=0$, hence $\dim(\ker\widetilde T_1)=1$.
Let us denote by $\tilde\phi\in\lW$ a function that forms a basis of $\ker\widetilde T_1$. 

If $\alpha$ does not have any zeros in the open interval $(a,b)$,
then $\tilde\phi$ is just a non-trivial solution of
$$
-(\alpha\tilde\phi)' + (\overline{\beta}+\alpha')\tilde\phi=0
$$
in $(a,b)$.

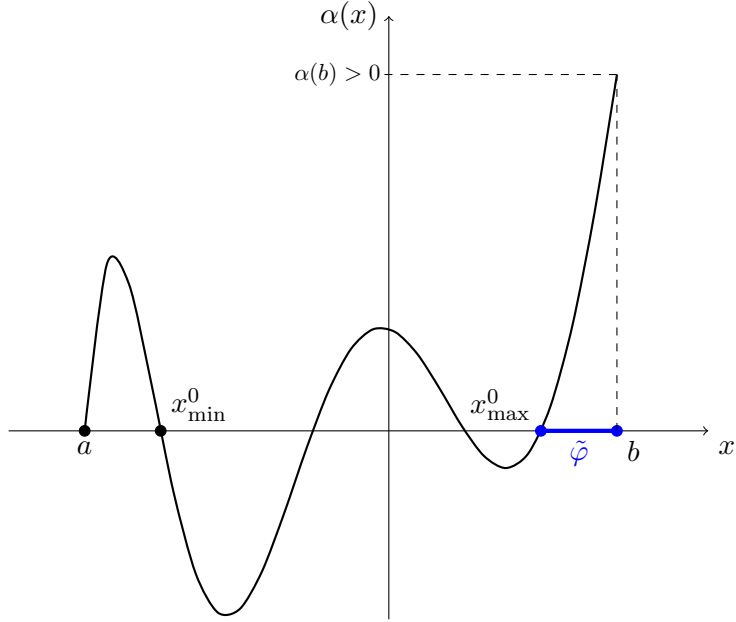
\begin{figure}[h]
\begin{tikzpicture}
  \draw[->] (-5, 0) -- (4.2, 0) node[below right] {$x$};
  \draw[->] (0, -2.5) -- (0, 5.5) node[left] {$\alpha(x)$};

  \draw[scale=1, thick, domain=-4:3, smooth, variable=\x]  plot ({\x}, {0.2*(\x/2+2)*(\x/2+1.5)*(\x/2+0.5)*(\x/2-0.5)*(\x/2-1)*(\x/2-3)^2});
 
  \coordinate (A) at (-4,0);
  \coordinate (B) at (3,0);
  \coordinate (Xmax) at (2,0);
  \coordinate (Xmin) at (-3,0);
  \coordinate (C) at (3,4.725);
 
  \filldraw (A) circle (2pt) node[below] {$a$};
  \filldraw[blue] (B) circle (2pt) node[below right, black] {$b$};
  \filldraw[blue] (Xmax) circle (2pt) 
    node[above left, black] {$x_\mathrm{max}^0$};
  \filldraw (Xmin) circle (2pt) node[above right] {$x_\mathrm{min}^0$};

  \draw[ultra thick, blue] (Xmax)--(B);
  \draw[dashed] (B)--(C);
  \draw[dashed] (C)--(-0.1,4.725);
  
  \node at (2.5,0) [below, blue] (phi) {$\tilde\phi$}; 

  \node at (0,4.725) [left, scale=0.8] (alpha) {$\alpha(b)>0$};
\end{tikzpicture}
\caption{For $\alpha$ satisfying $\alpha(a)=0$ and $\alpha(b)>0$ 
we denoted on the graph points $x_\mathrm{min}^0$ and
$x_\mathrm{max}^0$. The bold blue line segment contains the support
of $\tilde\phi$.}
\label{fig:case1}
\end{figure}

On the other hand, if $\alpha^{-1}(\{0\})\cap(a,b)\neq\emptyset$,
let us define
\begin{equation}\label{eq:xminxmax}
x_\mathrm{min}^0:=\min\Bigl(\alpha^{-1}(\{0\})\cap(a,b)\Bigr) \;,
    \quad x_\mathrm{max}^0:=\max\Bigl(\alpha^{-1}(\{0\})
    \cap(a,b)\Bigr) \;.
\end{equation}
Since in particular $\tilde\phi$ should satisfy 
the differential equation above in $(a,x_\mathrm{max}^0)$, 
where we have $\alpha(a)=\alpha(x_\mathrm{max}^0)=0$,
by the conclusion of the first subcase ($\alpha(a)=\alpha(b)=0$) 
we have that 
a.e.~$\tilde\phi|_{(a,x_\mathrm{max}^0)}=0$.
Thus, $\supp\tilde\phi\subseteq [x_\mathrm{max}^0,b]$
(see Figure \ref{fig:case1}).

In other cases the only differences are whether 
$\dim(\ker T_1)=1$ or $\dim(\ker\widetilde T_1)=1$,
and whether a function forming a basis is supported 
in $[a,x_\mathrm{min}^0]$ or $[x_\mathrm{max}^0,b]$.

\subsection{Case 2: $\alpha(a)\alpha(b)<0$}

In this case 
$\mathcal{W}_0=\{u\in \mathcal{W}:u(a)=u(b)=0\}$
(see Lemma \ref{lm:W0}).

Assume first that $\alpha(a)>0$ and $\alpha(b)<0$. 
Then for any $u\in\lW$ we have
\begin{equation*}
\dup{\lW'}{Du}{u}\lW = \alpha(b) |u(b)|^2 - \alpha(a)|u(a)|^2 \leq 0 
\,.
\end{equation*}
Hence, by Theorem \ref{tm:well-posedness} and with similar reasoning
as in the previous case we get 
that $(T_0,\widetilde T_1)=(T_1|_{\lW_0},\widetilde T_1|_\lW)$
is the only pair of mutually adjoint bijective realisations 
relative to $(T,\widetilde T)$.

Analogously, for $\alpha(a)<0$ and $\alpha(b)>0$
we have that $(T_1,\widetilde T_0)$ is the only pair of
mutually adjoint bijective realisations relative to 
$(T,\widetilde T)$.

Therefore, although in this case 
$\dim(\ker T_1)+\dim(\ker \widetilde T_1)=2$ 
(see Lemma \ref{lm:codim} and Remark \ref{rem:codim}(ii)),
we have only one bijective realisation. 
Hence, by Theorem \ref{tm:existence-ref-op}(ii)
and the analysis above, for $\alpha(a)>0$ we have
$(\dim(\ker T_1), \dim(\ker\widetilde T_1))=(2,0)$,
while for $\alpha(a)<0$ it is
$(\dim(\ker T_1), \dim(\ker \widetilde T_1))=(0,2)$.

Let us focus on the case $\alpha(a)<0$ and let us
study $\ker\widetilde T_1$. 
We define $x_\mathrm{min}^0$ and $x_\mathrm{max}^0$ 
as in \eqref{eq:xminxmax}.
They are well-defined since $\alpha(a)\alpha(b)<0$ and 
$\alpha$ is continuous, hence $\alpha^{-1}(\{0\})$ is not empty.
With the same argument as in the previous case we can conclude 
that for any $\tilde\phi\in\ker\widetilde T_1$ we have that
a.e.~$\tilde\phi|_{[x_\mathrm{min}^0,x_\mathrm{max}^0]}=0$.
Moreover, on both subintervals $(a,x_\mathrm{min}^0)$
and $(x_\mathrm{max}^0,b)$ we are in the same case 
regarding \eqref{eq:case1-VVtilda},
and this is precisely the reason why we have that one kernel 
is trivial, while the other being two-dimensional. 

\begin{figure}[ht]
\begin{tikzpicture}
  \draw[->] (-3, 0) -- (4.2, 0) node[below right] {$x$};
  \draw[->] (0, -2) -- (0, 5.5) node[left] {$\alpha(x)$};

  \draw[scale=1, thick, domain=-1.5:3, smooth, variable=\x]  plot ({\x}, {0.2*(\x/2+2)*(\x/2+1.5)*(\x/2+0.5)*(\x/2-0.5)*(\x/2-1)*(\x/2-3)^2});
 
  \coordinate (A) at (-1.5,0);
  \coordinate (B) at (3,0);
  \coordinate (Xmax) at (2,0);
  \coordinate (Xmin) at (-1,0);
  \coordinate (C) at (3,4.725);
  \coordinate (D) at (-1.5,-1.44196);
 
  \filldraw[red] (A) circle (2pt) node[below left, black] {$a$};
  \filldraw[blue] (B) circle (2pt) node[below right, black] {$b$};
  \filldraw[blue] (Xmax) circle (2pt) 
    node[above left, black] {$x_\mathrm{max}^0$};
  \filldraw[red] (Xmin) circle (2pt) node[below right, black] {$x_\mathrm{min}^0$};

  \draw[ultra thick, blue] (Xmax)--(B);
  \draw[ultra thick, red] (Xmin)--(A);
  \draw[dashed] (B)--(C);
  \draw[dashed] (C)--(-0.1,4.725);
  \draw[dashed] (A)--(D);
  \draw[dashed] (D)--(0.1,-1.44196);
  
  \node at (2.5,0) [below, blue] (phi) {$\tilde\phi_2$};
  \node at (-1.25,0) [above, red] (phi) {$\tilde\phi_1$};

 \node at (0,4.725) [left, scale=0.8] (alpha) {$\alpha(b)>0$};
 \node at (0,-1.44196) [right, scale=0.8] (alpha) {$\alpha(a)<0$};
\end{tikzpicture}
\caption{For $\alpha$ satisfying $\alpha(a)<0$ and $\alpha(b)>0$ 
we denoted on the graph points $x_\mathrm{min}^0$ and
$x_\mathrm{max}^0$. The bold red and blue line segments contain 
supports of $\tilde\phi_1$ and $\tilde\phi_2$, respectively.}
\label{fig:case2}
\end{figure}
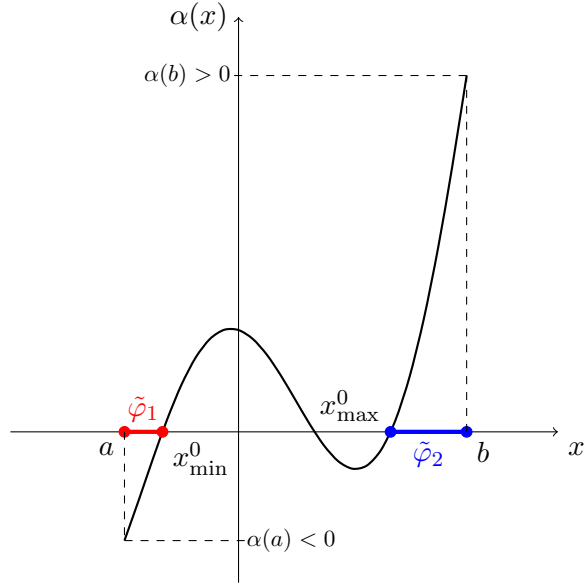

Thus, if we take $\tilde\phi_1,\tilde\phi_2 \in\lW$
such that $\tilde\phi_2=0$ on $[a,x_\mathrm{max}^0]$
and in $(x_\mathrm{max}^0,b)$ to be a non-trivial solution 
to the corresponding differential equation,
while $\tilde\phi_1=0$ on $[x_\mathrm{min}^0,b]$
and in $(a,x_\mathrm{min}^0)$ to be a non-trivial solution 
to the corresponding differential equation
(see Figure \ref{fig:case2}), 
then $\{\tilde\phi_1,\tilde\phi_2\}$ is a basis for
$\ker\widetilde T_1$.

\subsection{Case 3: $\alpha(a)\alpha(b)>0$}

As in the previous case, here we have
$\mathcal{W}_0=\{u\in \mathcal{W}:u(a)=u(b)=0\}$,
and the boundary operator reads (see \eqref{eq:D-scalar1d}):
\begin{equation*}
\dup{\lW'}{Du}{v}\lW = \alpha(b)u(b)\overline{v(b)}
    - \alpha(a)u(a)\overline{v(a)} \;, \quad
    u,v\in\lW \;.
\end{equation*}

Let us define
\begin{equation*}
\lV := \biggl\{ u\in\lW : u(b)= \sqrt{\frac{\alpha(a)}{\alpha(b)}}
    u(a) \biggr\} \;.
\end{equation*}
We are going to prove that the pair of subspaces $(\lV,\lV)$
satisfies condition (V1) and that $\lV=\lV^{[\perp]}$.
Then by Theorem \ref{tm:well-posedness} we shall have that 
operators $T_\mathrm{r}$ and $T_\mathrm{r}^*$, where 
$T_\mathrm{r}:=T_1|_{\lV}$, form a mutually adjoint pair of 
bijective realisations relative to 
$(T,\widetilde T)$. 

For an arbitrary $u\in\lV$ and $v\in\lW$ we have
\begin{equation}\label{eq:case3-D}
\begin{aligned}
\dup{\lW'}{Du}{v}\lW 
    &=\alpha(b)u(b)\overline{v(b)}-\alpha(a)u(a)\overline{v(a)}\\
&=\alpha(b)\biggl(u(b)\overline{v(b)}
    -\sqrt{\frac{\alpha(a)}{\alpha(b)}}u(a)
    \sqrt{\frac{\alpha(a)}{\alpha(b)}}\overline{v(a)}
    \biggr)\\
&=\alpha(b)u(b)\overline{\biggl(v(b)
    -\sqrt{\frac{\alpha(a)}{\alpha(b)}}
    v(a)\biggr)} \;.
\end{aligned}
\end{equation}
In particular,
\begin{equation*}
(\forall u,v\in\lV) \qquad 
\dup{\lW'}{Du}v\lW = 0 \,,
\end{equation*}
implying that $(\lV,\lV)$ satisfies condition (V1)
and that $\lV\subseteq \lV^{[\perp]}$.
Thus, it is left to show that $\lV^{[\perp]}\subseteq \lV$.

Let $v\in\lV^{[\perp]}$. Then by \eqref{eq:case3-D} for 
any $u\in\lV$ we have
$$
\alpha(b)u(b)\overline{\biggl(v(b)
    -\sqrt{\frac{\alpha(a)}{\alpha(b)}}
    v(a)\biggr)} = 0 \;.
$$
Since $\alpha(b)\neq 0$ and there exists $u\in\lV$
such that $u(b)\neq 0$
(e.g.~just consider a linear function),
this implies $v(b)=\sqrt{\frac{\alpha(a)}{\alpha(b)}}v(a)$, 
i.e.~$v\in\lV$.

Therefore, $(T_\mathrm{r},T_\mathrm{r}^*)$ is indeed 
a mutually adjoint pair of bijective realisations relative to 
$(T,\widetilde T)$.
It is evident that $\lW_0\subsetneqq\lV\subsetneqq\lW$,
hence by Theorem \ref{tm:existence-ref-op}(ii)
there are infinitely many bijective realisations. 
In particular, using the same theorem, 
we can conclude that both $\dim(\ker T_1)$
and $\dim(\ker\widetilde T_1)$ are greater or equal to 1. 
Now Remark \ref{rem:codim}(ii) implies that in fact we have
$\dim(\ker T_1) = \dim(\ker\widetilde T_1) = 1$.

Let $ \phi$ and $\tilde\phi$ span $\ker T_1$ and 
$\ker \widetilde T_1$, respectively.
We shall elaborate on explicit forms of functions
$\phi$ and $\tilde\phi$ later. 
First we want to determine all bijective realisations
relative to $(T,\widetilde T)$
(in terms of $\phi$ and $\tilde\phi$)
following Grubb's general extension theory,
which is concisely summarised in Appendix.
More precisely, we take $(T_\mathrm{r},T_\mathrm{r}^*)$
as the reference operators. 
In this part we shall mainly keep the notation used in
Appendix.
Let us remark that this procedure has already been done
in \cite[Section 6]{AEM-2017} for constant coefficients
$\alpha$ and $\beta$ (see also \cite{ABmn}). 

For any $u\in\lW$ there exist unique 
$u_\mathrm{r}\in\lV$ and $u_\mathrm{k}\in\ker T_1$
such that $u=u_\mathrm{r}+u_\mathrm{k}$
(see part (i) of Appendix). 
Moreover, $u_\mathrm{k}$ is of the form 
$C_u\phi$, so using 
\begin{align*}
u(a) &= u_\mathrm{r}(a) + C_u \phi(a) \\
u(b) &= u_\mathrm{r}(b) + C_u \phi(b)
\end{align*}
and $u_\mathrm{r}(b) = \sqrt{\frac{\alpha(a)}{\alpha(b)}} u_\mathrm{r}(a)$,
we get
\begin{equation}\label{eq:case3-cu}
C_u=\frac{u(b)-\sqrt{\frac{\alpha(a)}{\alpha(b)}}u(a)}
    {\phi(b)-\sqrt{\frac{\alpha(a)}{\alpha(b)}}\phi(a)}
\end{equation}
(note that the denominator is not equal to zero 
since $\phi\in\ker T_1\setminus\{0\}$ and 
$\ker T_1\cap \lV=\{0\}$).
Thus, the corresponding non-orthogonal projection
$p_\mathrm{k}:\lW\to\ker T_1$ is equal to
$p_\mathrm{k}(u)=C_u\phi$.
Similarly, $p_\mathrm{\tilde k}:\lW\to \ker\widetilde T_1$
is given by $p_\mathrm{\tilde k}(u)=\tilde C_u \tilde\phi$,
where
\begin{equation*}
\tilde C_u=\frac{u(b)-\sqrt{\frac{\alpha(a)}{\alpha(b)}}u(a)}
    {\tilde\phi(b)-\sqrt{\frac{\alpha(a)}{\alpha(b)}}\tilde\phi(a)} \;.
\end{equation*}

Since we seek for \emph{bijective} realisations, 
operator $B$ from part (ii) of Appendix 
should be bijective as well according to part (iii). 
Both kernels of $T_1$ and $\widetilde T_1$ are 
one-dimensional, hence the only (non-trivial) choice
is $\dom B = \lZ = \ker T_1$ and $\widetilde\lZ
=\dom{\widetilde T_1}$
(then also $\dom B^*=\ker\widetilde T_1$).
Then there exists $(c+id)\in\C$ such that 
$B\phi =(c+id)\tilde\phi$.
Therefore, all bijective realisations are indexed by 
$c+id\in\C\setminus\{0\}$
(for these values $B$ is an isomorphism).

The operator corresponding to $B$ we denote 
by $T_{c,d}=T_B$.
Recall that $T_0\subseteq T_{c,d}\subseteq T_1$.
By part (ii) of Appendix (see \eqref{eq:domAA*}),
$u\in\lW$ belongs to $\dom T_{c,d}$ if and only if
\begin{equation}\label{eq:case3-domTBa}
  P_{\ker\widetilde T_1}(T_1 u) = B(p_\mathrm{k}u) \,,  
\end{equation}
where $P_{\ker\widetilde T_1}$ is the orthogonal 
projection from $\lH$ onto $\ker \widetilde T_1$.

Let $u\in\lW$.
The right hand side of the equality above is equal 
to $(c+id)C_u \tilde\phi$, where $C_u$ is given by
\eqref{eq:case3-cu}.
For the left hand side we have
\begin{align*}
P_{\ker\widetilde T_1}(T_1 u) &= \frac{1}{\|\tilde\phi\|^2}
    \scp{T_1 u}{\tilde\phi} \tilde\phi \\
&= \frac{1}{\|\tilde\phi\|^2} \dup{\lW'}{Du}{\tilde\phi}\lW \tilde\phi \\
&= \frac{1}{\|\tilde\phi\|^2} 
    \Bigl(\alpha(b)u(b)\overline{\tilde{\phi}(b)}-\alpha(a)u(a)
    \overline{\tilde{\phi}(a)}\Bigr)\tilde\phi \;,
\end{align*}
where in the second equality we have used that
$\scp{u}{\widetilde T_1 \tilde\phi}=0$.
Thus, from \eqref{eq:case3-domTBa} we get that 
$u\in\lW$ belongs to $\dom T_{c,d}$ if and only if
\begin{equation}\label{eq:case3-domTBb}
\left( \frac{\alpha(b)\overline{\tilde\phi(b)}}{\|\tilde\phi\|^2}
    -\frac{(c+id)}{\phi(b)-\sqrt{\frac{\alpha(a)}{\alpha(b)}}
    \phi(a)} \right)u(b)
    =\left( \frac{\alpha(a)\overline{\tilde\phi(a)}}
    {\|\tilde\phi\|^2}-\frac{(c+id)\sqrt{\frac{\alpha(a)}
    {\alpha(b)}}}{\phi(b)-\sqrt{\frac{\alpha(a)}{\alpha(b)}}
    \phi(a)} \right) u(a) \;.
\end{equation}

Similarly using the second identity in \eqref{eq:domAA*}
(or directly computing the domain of the adjoint), 
we obtain that $u\in \mathcal{W}$ is in $\dom T^*_{c,d}$ if and only if 
\begin{equation}\label{eq:case3-domTB*}
\left( \alpha(b)\overline{\phi(b)}
    - \frac{\|\tilde\phi\|^2(c-id)}{\tilde\phi(b)-\sqrt{\frac{\alpha(a)}
    {\alpha(b)}}\tilde\phi(a)} \right) u(b)
    =\left( \alpha(a)\overline{\phi(a)}
    - \frac{\|\tilde\phi\|^2(c-id)\sqrt{\frac{\alpha(a)}{\alpha(b)}}}
    {\tilde\phi(b)-\sqrt{\frac{\alpha(a)}{\alpha(b)}}\tilde\phi(a)} \right) u(a) \;.
\end{equation} 

Therefore, the set of all pairs of mutually adjoint bijective realisations 
relative to $(T,\widetilde T)$ in this case is given by
\begin{equation}\label{eq:case3-all}
\Bigl\{(T_{c,d},T_{c,d}^*) : c,d\in\R^2\setminus\{(0,0)\}\Bigr\} \bigcup 
    \,\bigl\{(T_\mathrm{r},T_\mathrm{r}^*)\bigr\} \;.
\end{equation}
All bijective realisations are parametrised by \emph{one}
complex parameter ($c+id$), which is in parallel to the fact that 
the dimension of both kernels $\ker T_1$ and $\ker \widetilde T_1$
is \emph{one}.

Note that $\dom T_{c,d}=\lW_0+\ker\widetilde T_1$ (see
Corollary \ref{cor:decomposition}) if and only if 
$\tilde\phi\in\dom T_{c,d}$. Indeed, then 
$\lW_0+\ker\widetilde T_1\subseteq \dom T_{c,d}$ and the inclusion cannot be
strict as in that case it would be impossible that both operators 
$T_{c,d}$ and $T_1|_{\lW_0+\ker\widetilde T_1}$ are bijective.
From the above it can be easily seen that $\tilde\phi\in\dom T_{c,d}$
is achieved if and only if 
\begin{equation*}
c+id = \frac{\dup{\lW'}{D\tilde\phi}{\tilde\phi}\lW}{\|\tilde\phi\|^2 C_{\tilde\phi}} \;.
\end{equation*}

Let us go back to kernels of $T_1$ and $\widetilde T_1$, so that we can 
derive some properties of functions $\phi$ and $\tilde\phi$.

If $\min_{x\in[a,b]}|\alpha(x)|>0$, then we get $\phi$ and 
$\tilde\phi$ simply by taking non-trivial solutions of
\begin{equation}\label{eq:case3-kernels}
(\alpha \phi)' + \beta\phi = 0 \qquad\hbox{and}\qquad
-(\alpha\tilde\phi)' + (\overline{\beta}+\alpha')\tilde\phi=0
\end{equation}
on $(a,b)$. Thus, a possible choice is ($x\in [a,b]$):
\begin{equation}\label{eq:case3-kernels-basis}
\phi(x) = \frac{1}{\alpha(x)}\operatorname{exp}
    \Bigl(-\int\frac{\beta(x)}{\alpha(x)}\,dx\Bigr)
    \qquad \hbox{and}\qquad
    \tilde\phi(x) = \operatorname{exp}\Bigl(
    \int\frac{\overline{\beta(x)}}{\alpha(x)}\,dx\Bigr) \;.
\end{equation}

If $\alpha^{-1}(\{0\})\cap(a,b)$ is not empty, we define
$x_\mathrm{min}^0$ and $x_\mathrm{max}^0$ as in 
\eqref{eq:xminxmax}.
Here we can apply the same inference as in Case 1 to conclude
that functions $\phi$ and $\tilde\phi$ are supported 
on $[a,x_\mathrm{min}^0]$ or $[x_\mathrm{max}^0,b]$,
while on the supports we just use \eqref{eq:case3-kernels-basis}
(one needs to be aware that now integrals are improper, but for sure convergent
as we know that such non-trivial $\phi$ and $\tilde\phi$
should exist in $\lW$).

\begin{figure}[H]
\begin{tikzpicture}
  \draw[->] (-5, 0) -- (4.2, 0) node[below right] {$x$};
  \draw[->] (0, -2.5) -- (0, 5.5) node[left] {$\alpha(x)$};

  \draw[scale=1, thick, domain=-3.8:3, smooth, variable=\x]  plot ({\x}, {0.2*(\x/2+2)*(\x/2+1.5)*(\x/2+0.5)*(\x/2-0.5)*(\x/2-1)*(\x/2-3)^2});
 
  \coordinate (A) at (-3.8,0);
  \coordinate (B) at (3,0);
  \coordinate (Xmax) at (2,0);
  \coordinate (Xmin) at (-3,0);
  \coordinate (C) at (3,4.725);
  \coordinate (D) at (-3.8,1.87163);
 
  \filldraw[red] (A) circle (2pt) node[below left, black] {$a$};
  \filldraw[blue] (B) circle (2pt) node[below right, black] {$b$};
  \filldraw[blue] (Xmax) circle (2pt) 
    node[above left, black] {$x_\mathrm{max}^0$};
  \filldraw[red] (Xmin) circle (2pt) node[above right, black] {$x_\mathrm{min}^0$};

  \draw[ultra thick, blue] (Xmax)--(B);
  \draw[ultra thick, red] (Xmin)--(A);
  \draw[dashed] (B)--(C);
  \draw[dashed] (C)--(-0.1,4.725);
  \draw[dashed] (A)--(D);
  \draw[dashed] (D)--(0.06,1.87163);
  
  \node at (2.5,0) [below, blue] (phi) {$\tilde\phi$};
  \node at (-3.4,0) [below, red] (phii) {$\phi$}; 

  \node at (0,4.725) [left, scale=0.8] (alpha) {$\alpha(b)>0$};
  \node at (0,1.87163) [right, scale=0.8] (alpha) {$\alpha(a)>0$};
\end{tikzpicture}
\caption{For $\alpha$ satisfying $\alpha(a)>0$ and $\alpha(b)>0$ 
we denoted on the graph points $x_\mathrm{min}^0$ and
$x_\mathrm{max}^0$. The bold red and blue line segments contain 
supports of $\phi$ and $\tilde\phi$, respectively.}
\label{fig:case3}
\end{figure}
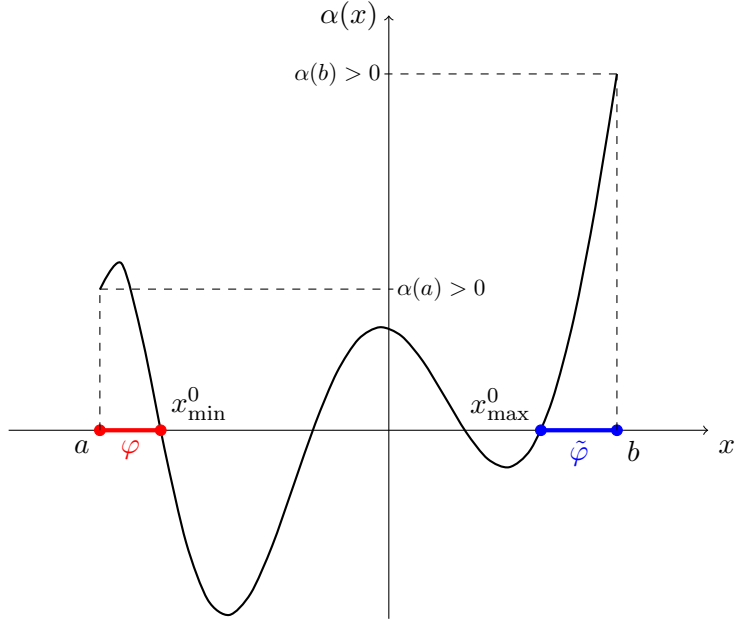

To be more specific, let us assume that $\alpha(a)>0$ and $\alpha(b)>0$.
Then any solution in $\lW$ of the first equation in 
\eqref{eq:case3-kernels} must satisfy 
$\phi|_{[x_\mathrm{min}^0,b]}=0$, while for the second equation 
we have $\tilde\phi|_{[a,x_\mathrm{max}^0]}=0$ 
(see Figure \ref{fig:case3}). 
In particular, under this assumption we have $\phi(b)=\tilde\phi(a)=0$, 
which could be used to simplify \eqref{eq:case3-domTBb}
and \eqref{eq:case3-domTB*}.
This also implies that $\lW_0+\ker\widetilde T_1 = \{u\in\lW : u(a)=0\}$
and $\lW_0+\ker T_1 = \{u\in\lW : u(b)=0\}$.

\begin{remark}
    Equation (5.7) covers all (linear) boundary conditions of the form $\gamma u(b)=\delta u(a)$, where $(\gamma ,\delta) \in \C^2\setminus \{(0,0)\}$, 
except the one that is satisfied by all functions from $\ker T_1$ (and then also $\phi$).
To justify this claim let us just study the case $c=d=0$ (the only case which does not lead to a bijective realisation).
We get
\begin{align*}
    \alpha (b)\overline{\tilde\phi(b)}u(b)=\alpha(a)\overline{\tilde\phi(a)}u(a)\;,
\end{align*}
implying $\iscp{u}{\tilde \phi}=0$,
which concludes to $u\in \lW_0\dot{+}\ker T_1$ using \eqref{eq:decomposition}. Thus, the above boundary condition is satisfied by functions from $\ker T_1$. The proof of the remaining part, that is the fact that all other boundary conditions are attained, is left to the readers.

Using Lemma \ref{lm:con}, all bijective realisations can be characterised in a more concise way. Indeed, all possible domains of the bijective realisations are given by
\begin{align*}
    \lV=\lW_0\,\dot{+}\,\span\{\tilde \phi +\lambda \phi\} \;,\quad \lambda \in \C\;.
\end{align*}
This can serve as another evidence to the above claim.

However, the approach using the universal classification theory has some additional advantages 
when studying e.g.~the spectrum and the resolvent of realisations. Moreover, once the classification is established, 
choosing the desired properties for realisations comes down to choosing the same properties for operator $B$,
which is often easier to control.
\end{remark}

\begin{remark}
Instead of Grubb's general extension theory, one could use the boundary triplets formalism 
\cite[Section 13.4]{GG}, which could be even more convenient for the calculation of resolvents.
\end{remark}

\begin{remark}
    In \cite[Section 6.1]{AEM-2017} an analysis 
    of operators \eqref{eq:TTtilda} is made for
    $\alpha=\beta=1$, which is a special case of the situation 
    described in Case 3 above.
    Besides classifying all bijective realisations, in the aforementioned
    reference a further step is made to distinguish \emph{bijective realisations
    with signed boundary map}, i.e.~those which are provided by 
    Theorem \ref{tm:well-posedness}.
    The same can be conducted here using \cite[Theorem 18(ii)]{AEM-2017}.
    Here we have decided to omit that part as already bijective realisations
    provide well-posedness of the corresponding abstract problems.
\end{remark}

\subsection{Summary}
   
Depending on the values of $\alpha$ at end-points,
the pairs of subspaces $(\mathcal{V},\widetilde{\mathcal{V}})$
for which we obtain bijective realisations,
i.e.~such that $(T_1|_\lV, \widetilde T_1|_{\widetilde \lV})$
is a pair of mutually adjoint bijective realisations relative to 
$(T,\widetilde T)$, where $T$ and $\widetilde T$
are given by \eqref{eq:TTtilda}, are:

\begin{align*}
    \begin{tabular}{|c|c|c|c|}
  \hline
  $\alpha$ at end-points & No.~of bij.~realisations & $(\mathcal{V},\widetilde{\mathcal{V}})$ \\
  \hline   
   $\alpha(a)\alpha(b)\leq 0$ & $1 $ &
    \begin{tabular}{c|c} 
        $\alpha(a)\geq 0 \wedge \alpha(b)\leq 0$ & $(\lW_0,\lW)$\\
        \hline 
        $\alpha(a)\leq 0 \wedge \alpha(b)\geq 0$ & $(\lW,\lW_0)$  \\
        \end{tabular}\\
  \hline 
  $\alpha(a)\alpha(b)>0$ & $\infty$ & \eqref{eq:case3-all} (see \eqref{eq:case3-domTBb}--\eqref{eq:case3-domTB*}) \\
  \hline 
\end{tabular}
\end{align*}
Thus, a classification of bijective realisations is needed only in the case when 
$\alpha$ has the same sign in both end-points.

\section{Examples}\label{sec:examples}

\begin{enumerate}
    \item Take the interval $(0,2)$ and coefficients $\alpha(x)=1-x$
    and $\beta=1$. Then 
          \[T\phi=((1-x)\phi)'+\phi \]and 
          \[\widetilde{T}\phi = -((1-x)\phi)'\;.\]
          Here $2\Re \beta+\alpha'=1> 0$ on $(0,2)$,
          meaning that $(T,\widetilde T)$ is a pair of abstract Friedrichs 
          operators. Moreover, this example belongs to Case 2 of
          the previous section. 
          
          Using \eqref{eq:case3-kernels-basis} on $(0,1)$ and $(1,2)$ 
          separately we get that for $\phi\in\ker T_1$ necessarily 
\begin{align*}
    \phi=\left\{\begin{array}{lcl}
                       c_1 &,& \hbox{in} \  (0,1)  \\ 
                        c_2 &, & \hbox{in} \  (1,2)\;,
         \end{array}\right.
\end{align*}
for some constants $c_1,c_2\in\C$.
We have $\phi\in\lW$.
Indeed, it is evident that $\phi\in {L}^2(0,2)$, 
while for $\psi\in {C}^\infty_c(0,2)$ we have
\begin{align*}
     \int_0^2(1-x)\phi(x)\psi'(x)\,dx 
    & = c_1 \int_0^1(1-x)\psi'(x)\,dx+c_2 \int_1^2(1-x) \psi'(x)\,dx \\
    & = c_1\int_0^1 \psi(x) \,dx  +c_2 \int_1^2 \psi(x) \,dx\\
    & = \int_0^2 \phi(x) \psi(x)\, dx\;.
\end{align*}
This means $((1-x)\phi)'=-\phi \in {L}^2(0,2)$, thus $\phi\in \mathcal{W}$.
Therefore, $\dim\ker T_1=2$ (since we have two parameters in the definition
of $\phi$).

On the other hand, $\tilde\phi\in\ker\widetilde T_1$ implies      
\begin{align*}
    \tilde\phi(x)=\left\{\begin{array}{lcl}
                       \frac{d_1}{1-x} &,& x\in (0,1)  \\ 
                        \frac{d_2}{1-x} &,& x\in (1,2)\;,
         \end{array}\right.
\end{align*}
for some constants $d_1,d_2\in\C$. But it can easily be seen that 
$\tilde\phi\in {L}^2(0,2)$ if and only if $d_1=d_2=0$.
Hence, $\ker \widetilde{T}_1=\{0\}$ and $\dim \ker \widetilde{T}_1=0$,
justifying the results obtained in Case 2 of the previous section.

It is interesting to note that for $c_1\neq c_2$ we have
$\phi'\notin {L}^2(0,2)$, 
because $\phi'=(c_2-c_1)\delta_1$ 
(here $\delta_1$ is the Dirac measure at $1$) and so 
$\phi\notin {H}^1(0,2)$.
Thus, ${H}^1(0,2)\subsetneqq \mathcal{W}$.

Moreover, it is evident that $\tilde\phi\in {H}^1_{\loc}([0,2]\setminus\{1\})$
for any choice  of parameters $d_1, d_2$.
Indeed, for any  subinterval $[c,d]\subseteq [0,2]\setminus\{1\}$ 
we have $\tilde\phi|_{(c,d)}\in \mathrm{H}^1{(c,d)}$. 
Since $\tilde\phi\not\in\lW$ this shows that $\mathcal{W}$ is 
indeed a proper subspace of ${H}^1_{\loc}([0,2]\setminus\{1\})$,
i.e.~$\lW\subsetneqq {H}^1_{\loc}([0,2]\setminus\{1\})$.\\

    \item Take the same example as above, but now on the interval $(0,1).$ Here  using \eqref{eq:case3-kernels-basis} again we get 
that $\phi \in \ker T_1$ implies that $\phi= c$, for some 
constant $c\in\C$. 
Since $\phi\in H^1(0,1)$, it is contained in the graph space $\lW$. 
Hence $\ker T_1 =\operatorname{span}\{1\}$ and $\dim\ker T_1=1 $.

Furthermore, for $\tilde{\phi}\in \ker \widetilde{T}_1$ necessarily
 \begin{align*}
    \tilde{\phi}(x)=\frac{d}{(x-1)} \;, \quad x\in (0,1)\;, 
\end{align*}           
for some constant $d\in \mathbb{C}$. 
But $\tilde{\phi}\in {L}^2(0,1)$ if and only if $d=0$.
Hence, $\ker \widetilde{T}_1=\{0\}$ and $\dim \ker \widetilde{T}_1=0$.

This coincides with the results obtained in Case 1 of the previous section.\\

\item Let us consider another example that fits into the setting of Case 1 of the previous section.
Take $\alpha(x)=x(x-1)$ and $\beta = 1$ on the interval $(0,1)$.
Here $\alpha'(x)=2x-1$, so we have $2\Re{\beta}+\alpha'\geq 1>0$ in $(0,1)$. By \eqref{eq:case3-kernels-basis},
$\phi \in \ker T_1$ and $\tilde \phi \in \ker \widetilde{T}_1$ imply
   \[\phi(x) = \frac{c}{(x-1)^2} \;,\qquad \tilde \phi(x) = d\left (\frac{x-1}{x}\right ) \;, \]
   for some constants $c,d \in \mathbb{C}$. But $\phi,\tilde{\phi} \in {L}^2(0,1)$ if and only if $c=d=0$.
   Hence, $\ker T_1 = \ker \widetilde T_1=\{0\}$.\\ 

\item Take $\alpha(x) = (x-1)(x-2)$ and $\beta=2$ on the interval $(0,3)$. Then $\alpha $ has two zeroes on the interval $(0,3)$. 
Here $\alpha'(x)=2x-3$, hence we have $2\Re{\beta}+\alpha'\geq 1>0$ in $(0,3)$. Again using \eqref{eq:case3-kernels-basis} on subintervals
$(0,1)$, $(1,2)$ and $(2,3)$ separately we get that 
$\phi\in \ker{T}_1$ implies 
\begin{align*}
    \phi(x)=
    \left\{\begin{array}{lcl}
        c_1\left(\frac{x-1}{x-2}\right)^2  &,& x\in(0,1) \\
        c_2 \left(\frac{x-1}{x-2}\right)^2 &,& x\in(1,2) \\
        c_2 \left(\frac{x-1}{x-2}\right)^2 &,& x\in (2,3) \;,
        \end{array}\right.
\end{align*}
for some constants $c_1,c_2,c_3\in \mathbb{C}$.
But $\phi\in {L}^2(0,3)$ if and only if $c_2=c_3=0$. 
Moreover, for $c_2=c_3=0$ we have $\phi \in \mathcal{W}$,
implying $\dim \ker T_1 =1$. 

On the other hand, $\tilde \phi \in \ker \widetilde{T}_1$ implies
\begin{align*}
    \tilde{\phi}(x)=
    \left\{\begin{array}{lcl}
        d_1 \left(\frac{x-2}{x-1}\right)^2  &,& x\in(0,1)\\
        d_2 \left(\frac{x-2}{x-1}\right)^2  &,& x\in(1,2)\\
        d_3 \left(\frac{x-2}{x-1}\right)^2  &,& x\in (2,3) \;,
        \end{array}\right. 
\end{align*}
for some constants $d_1,d_2,d_3\in \mathbb{C}$. But $\tilde \phi \in {L}^2(0,3)$ if and only if $d_1=d_2=0$, and for $d_1=d_2=0$ we have 
$ \tilde{\phi} \in \mathcal{W}$. So, $\dim\ker \widetilde{T}_1=1$,
which is in accordance with Case 3 of the previous section.

\end{enumerate}

\section{Concluding remarks}\label{sec:conclusion}

In Section \ref{sec:1d-class} we provided a full classification 
of \emph{all} bijective realisations relative to \eqref{eq:TTtilda}.
Compared to \cite[Section 6.1]{AEM-2017}, here we were able to treat
operators with variable coefficients, thus indeed cover scalar ($r=1$)
one-dimensional ($d=1$) classical Friedrichs operators (CFO)
in the full generality.
One needs to be aware that by \emph{all}\/ we mean among all 
possible \emph{linear} boundary conditions, as this whole
theory is linear. Of course, by a standard procedure one can easily
generalise these results to include
inhomogeneous (linear) boundary conditions as well
(for some examples of classical Friedrichs operators a 
well-posed variational theory has been developed recently \cite{BH21}
allowing directly inhomogeneous boundary conditions).

The decomposition \eqref{eq:decomposition}, developed in Section 
\ref{sec:decomposition}, can be seen as a von Neumann type formula for non-symmetric operators.
It has been used in the analysis of (CFO)
to simplify some arguments and to make all the results more convincing, 
but it was not essential.
However, we are certain that in the study of (CFO) (or some other operators) in the vectorial 
and/or higher-dimensional setting it will play more important role
(maybe even fundamental).
Indeed, in such a generality it is hard to expect that all objects will be explicit as
in Section \ref{sec:1d-class} (e.g.~in the higher-dimensional case we expect kernels
$\ker T_1$ and $\ker \widetilde{T}_1$ to be inifinite-dimensional), 
so general results as the decomposition \eqref{eq:decomposition}
should be of a big help.

Another perspective of decomposition \eqref{eq:decomposition}
might be in better understanding of three equivalent descriptions
for boundary conditions of abstract Friedrichs operators \cite{EGC}.
More precisely, equivalence between all abstract descriptions was closed
in \cite{ABcpde}, but their explicit relation was obtained only in special situations
when certain projectors exist (see also \cite{ABjde, ABVisrn}).
These relations are important as e.g.~\emph{strongly enforced} boundary conditions
(those which are incorporated in the definition of the solution space) are not 
always desirable (e.g.~in certain numerical schemes). Therefore, 
an explicit method of switching from one formulation to another is preferable.
We hope that with \eqref{eq:decomposition} at hand these results can be improved.

\section{Appendix: Grubb's classification}

In this section we briefly recall the general extension
theory of (closed and) densely defined operators on Hilbert
spaces following \cite[Chapter 13]{GG}
(see also \cite{Grubb-1968} where the main results were
already obtained). 

Let $(A_0,\widetilde{A}_0)$ and $(A_1,\widetilde{A}_1)$ be 
two pairs of mutually adjoint, closed and densely defined operators 
on $\lH$ satisfying
\begin{equation*}
	A_0\;\subseteq\;(\widetilde{A}_0)^*\;=\;A_1
	\qquad\textrm{and}\qquad \widetilde{A}_0\;\subseteq\;(A_0)^*\;=\;\widetilde{A}_1\,,
\end{equation*}
which admit a further pair $(A_\mathrm{r},A_\mathrm{r}^*)$ of reference operators that are closed, satisfy 
$A_0\subseteq A_\mathrm{r}\subseteq A_1$,
equivalently $\widetilde{A}_0\subseteq A_\mathrm{r}^*
\subseteq \widetilde{A}_1$, and are invertible with everywhere defined bounded inverses $A_\mathrm{r}^{-1}$ and $(A_\mathrm{r}^*)^{-1}$.
Then the following holds.

	\begin{itemize}
		\item[(i)] There are decompositions 
		\begin{equation*}
		\dom A_1\;=\;\dom A_\mathrm{r} \dotplus\ker A_1\qquad \hbox{and} \qquad  \dom \widetilde{A}_1\;=\;\dom A_\mathrm{r}^*\dotplus\ker \widetilde{A}_1\,,
		\end{equation*}
		the corresponding (non-orthogonal) projections 
		\begin{equation*}
		\begin{array}{ll}
		p_\mathrm{r}:\dom A_1\to\dom A_\mathrm{r}\,, &\qquad  p_{\tilde{\mathrm{r}}}:\dom \widetilde{A}_1\to\dom A_\mathrm{r}^*\,,\\
		p_{\mathrm{k}}:\dom A_1\to\ker A_1\,, &\qquad p_{\tilde{\mathrm{k}}}:\dom \widetilde{A}_1\to\ker \widetilde{A}_1\,,
		\end{array}
		\end{equation*}
		satisfying 
		\begin{equation*}
		\begin{array}{ll}
		p_\mathrm{r}\;=\;A_\mathrm{r}^{-1}A_1\,, &\qquad p_{\tilde{\mathrm{r}}}\;=\;(A_\mathrm{r}^*)^{-1}
		    \widetilde{A}_1\,, \\
		p_{\mathrm{k}}\;=\;\mathbbm{1}-p_\mathrm{r}\,, &\qquad p_{\tilde{\mathrm{k}}}\;=\;\mathbbm{1}-p_{\tilde{\mathrm{r}}}\,,
		\end{array}
		\end{equation*}
		and being continuous with respect to the graph norms.
		\item[(ii)] There is a one-to-one correspondence between all pairs of mutually adjoint operators $(A,A^*)$ with $A_0\subseteq A\subseteq A_1$,
		equivalently $\widetilde{A}_0\subseteq A^*\subseteq \widetilde{A}_1$, and all pairs of densely defined mutually adjoint operators $B:\mathcal{Z}\to\widetilde{\mathcal{Z}}$ and
		$B^*:\widetilde{\mathcal{Z}}\to \mathcal{Z}$, with domains $\dom B\subseteq\mathcal{Z}$ and $\dom B^\ast\subseteq\widetilde{\mathcal{Z}}$,
		where $\mathcal{Z}$ and $\widetilde{\mathcal{Z}}$ run through all closed subspaces of $\ker A_1$ and $\ker \widetilde A_1$.
		The correspondence is given by
		\begin{equation}\label{eq:domAA*}
		\begin{split}
		\dom A\;&=\;\Bigl\{ u\in\dom A_1\, : \,p_\mathrm{k}u\in\dom B\,, \; P_{\widetilde{\mathcal{Z}}}(A_1u)=B (p_\mathrm{k}u) \Bigr\} \,,\\
		\dom A^*\;&=\;\Bigl\{v\in\dom \widetilde{A}_1\, : \,p_{\tilde{\mathrm{k}}}v\in\dom B^*\,, \; P_\mathcal{Z}(\widetilde{A}_1 v)=B^* 	
		(p_{\tilde{\mathrm{k}}}v) \Bigr\}\,,
		\end{split}
		\end{equation}
		and conversely, by
		\begin{equation*}
		\begin{array}{rclrclrcl}
		\dom B &=& p_\mathrm{k}\dom A \,, \quad &
		    \mathcal{Z} &=& \overline{\dom B}\,, \quad &
		    B (p_\mathrm{k}u) &=& P_{\widetilde{\mathcal{Z}}}(A_1u)\,, \\
		\dom B^* &=& p_{\tilde{\mathrm{k}}}\dom A^*\,, \quad &
		    \widetilde{\mathcal{Z}} &=& \overline{\dom B^*}\,, \quad &
		    B^* (p_{\tilde{\mathrm{k}}}v) &=& P_\mathcal{Z}(\widetilde{A}_1 v)\,,
		\end{array}
		\end{equation*}
		where $P_\mathcal{Z}$ and $P_{\widetilde{\mathcal{Z}}}$ are the \/\emph{orthogonal}\/ projections from $\lH$ onto $\mathcal{Z}$ and $\widetilde{\mathcal{Z}}$.
		\item[(iii)] In the correspondence above, $A$ is injective, resp.~surjective, resp.~bijective, if and only if so is $B$.
		\item[(iv)] When $A_B$ corresponds to $B$ as above, then
		\begin{equation*}
		\begin{split}
		\dom A_B\;&=\;\left\{w_0+(A_\mathrm{r})^{-1}(B\nu+\tilde\nu)+\nu\,\left|\!\!
		\begin{array}{c}
		w_0\in\dom A_0 \\
		\nu\in\dom B \\
		\tilde\nu\in\ker \widetilde{A}_1\ominus\widetilde{\mathcal{Z}}
		\end{array}\!\!
		\right.\right\} \,,\\
		A_B\bigl(w_0&+(A_\mathrm{r})^{-1}(B\nu+\tilde\nu)+\nu\bigr)\;=\;A_0 w_0+B\nu+\tilde\nu
		\end{split}
		\end{equation*}
		and
		\begin{equation*}
		\begin{split}
		\dom (A_B)^* \;&=\;\left\{\tilde{w}_0+(A_{\mathrm{r}}^*)^{-1}(B^*\tilde\mu+\mu)+\tilde\mu\,\left|\!\!
		\begin{array}{c}
		\tilde{w}_0\in\dom \widetilde{A}_0 \\
		\tilde\mu\in\dom B^* \\
		\mu\in\ker A_1\ominus\mathcal{Z}
		\end{array}\!\!
		\right.\right\} \,,\\
		(A_B)^*\bigl(\tilde{w}_0&+(A_{\mathrm{r}}^*)^{-1}(B^*\tilde\mu+\mu)+\tilde\mu\bigr)\;=\;\widetilde{A}_0 \tilde{w}_0+B^*\tilde\mu+\mu\,,
		\end{split}
		\end{equation*}
		and, moreover,
		\begin{equation*}
		(A_B)^*\;=\;A_{B^*}\,.
		\end{equation*}
		For the trivial choice $\mathcal{Z}=\widetilde{\mathcal{Z}}=\{0\}$ one has $A_B=A_\mathrm{r}$.
	\end{itemize}


\section{Acknowledgements}
We warmly thank N.~Antonić and A.~Michelangeli
for enlightening discussions on
the subject.

This work is supported by the Croatian Science Foundation under project 
IP-2018-01-2449 MiTPDE.

\end{document}